\documentclass[10pt]{book}
\usepackage{array}
\usepackage{color}
\usepackage{amssymb}
\usepackage{amsmath}
\usepackage{epsfig}
\usepackage[authoryear]{natbib}
\usepackage{setspace}
\usepackage{amsbsy}
\usepackage{multirow}

\newtheorem{theorem}{Theorem}[section]
\newtheorem{lemma}{Lemma}[section]
\newtheorem{remark}{Remark}[section]
\newtheorem{example}{Example}[section]

\begin{document}

\setcounter{page}{1}
%\begin{document}
%\pagestyle{fancy}
\renewcommand{\baselinestretch}{1.2}
%\lhead[\fancyplain{} \leftmark]{}
%\chead[]{}
%\rhead[]{\fancyplain{}\rightmark}
%\cfoot{}
%\headrulewidth=0pt
\markright{
%\hbox{\footnotesize\rm Statistica Sinica
%{\footnotesize\bf ??}(200?), 000-000}\hfill
}
\markboth{\hfill{\footnotesize\rm JIE YANG, ABHYUDAY MANDAL AND DIBYEN MAJUMDAR
}\hfill}
{\hfill {\footnotesize\rm OPTIMAL DESIGN FOR BINARY RESPONSE} \hfill}
\renewcommand{\thefootnote}{}
$\ $\par
\fontsize{10.95}{14pt plus.8pt minus .6pt}\selectfont
\vspace{0.8pc}
\centerline{\large\bf OPTIMAL DESIGNS FOR $2^k$ FACTORIAL EXPERIMENTS}
\vspace{2pt}
\centerline{\large\bf WITH BINARY RESPONSE}
\vspace{.4cm}
\centerline{Jie Yang$^{1}$, Abhyuday Mandal$^{2}$ and Dibyen Majumdar$^{1}$}
\vspace{.4cm}
\centerline{\it  $^1$University of Illinois at Chicago and $^2$University of Georgia}
\vspace{.55cm}
\fontsize{9}{11.5pt plus.8pt minus .6pt}\selectfont

\begin{quotation}
\begin{center} {\bf{ \it Abstract:}}\end{center}
   We consider the problem of obtaining D-optimal designs for factorial experiments with a binary response and $k$ qualitative factors each at two levels. We obtain a characterization for a design to be locally D-optimal. Based on this characterization, we develop efficient numerical techniques to search for locally D-optimal designs. Using prior distributions on the parameters, we investigate EW D-optimal designs, which are designs that maximize the determinant of the expected information matrix. It turns out that these designs can be obtained very easily using our algorithm for locally D-optimal designs and are very good surrogates for Bayes D-optimal designs. We also investigate the properties of fractional factorial designs and study the robustness with respect to the assumed parameter values of locally D-optimal designs.

\par

\vspace{9pt}  {\it Key words and phrases:} Generalized linear model, full factorial design, fractional factorial design, D-optimality, uniform design, EW D-optimal design.
\par
\end{quotation}\par

\fontsize{10.95}{14pt plus.8pt minus .6pt}\selectfont
\setcounter{chapter}{1}
\setcounter{equation}{0} %-1
 {\bf 1. Introduction}

Our goal is to determine optimal and efficient designs for factorial experiments with qualitative factors and a binary response. The traditional factorial design literature deals with experiments where the factors have discrete levels and the response follows a linear model (see, for example, Xu et al. (2009) and references therein). On the other hand, there is a growing body of literature on optimal designs for quantitative factors with binary or categorical response. For the specific experiments we study, however, the design literature is meager. Consequently, these experiments are usually designed by the guidelines of traditional factorial design theory for linear models. As we shall see, the resulting designs can be quite inefficient, especially when compared to designs that make use of prior information when it is available. Our goal is to address this problem directly and determine efficient designs specifically for experiments with qualitative factors and a binary response.

We assume that the process under study is adequately described by a generalized linear model (GLM). GLMs have been widely used for modeling binary response. Stufken and Yang (2012) noted that ``the study of optimal designs for experiments that plan to use a GLM is however not nearly as well developed (see also Khuri, Mukherjee, Sinha and Ghosh, 2006), and tends to be much more difficult than the corresponding and better studied problem for the special case of linear models.'' For optimal designs under GLMs, there are four different approaches proposed in the literature to handle the dependence of the design optimality criterion on the unknown parameters, (1) local optimality approach of Chernoff (1953) in which the parameters are replaced by assumed values; (2) Bayesian approach (Chaloner and Verdinelli (1995)) that incorporates prior belief on unknown parameters; (3) maximin approach that maximizes the minimum efficiency over a range of values of the unknown parameters (see Pronzato and Walter (1988) and Imhof (2001)); and (4) sequential approach where the design and parameter estimates are updated in an iterative way (see Ford, Titterington and Kitsos (1989)). In this paper, we will focus on local optimality and study D-optimal factorial designs under GLMs. We also consider Bayes optimality and study a surrogate for Bayes D-optimal designs that has many desirable properties.

The methods for analyzing data from GLMs have been discussed in depth in the literature (for example, McCullagh and Nelder (1989), Agresti (2002), Lindsey (1997), McCulloch and Searle (2001), Dobson and Barnett (2008) and Myers, Montgomery and Vining (2002)). Khuri, Mukherjee, Sinha and Ghosh (2006) provided a systematic study of the optimal design problem in the GLM setup and recently there has been an upsurge in research in both theory and computation of optimal designs. Russell et al.~(2009), Li and Majumdar (2008, 2009), Yang and Stufken (2009), Yang et al.~(2011), Stufken and Yang (2012) are some of the papers that developed theory and Woods et al.~(2006), Dror and Steinberg (2006, 2008), Waterhouse et al.~(2008), Woods and van de Ven (2011) focused on developing efficient numerical techniques for obtaining optimal designs under generalized linear models. Our focus is on optimal designs for GLMs with qualitative factors.

The special case of $2^2$ experiments with qualitative factors and a binary response was studied by Yang, Mandal and Majumdar (2012), where we obtained optimal designs analytically in special cases and demonstrated how to obtain a solution in the general case using cylindrical algebraic decomposition. The optimal allocations were shown to be robust to the choice of the assumed values of the model parameters. Gra{\ss}hoff and Schwabe (2008) has some relevant results for the $k=2$ factor case. The extension for $k>2$ factors is substantial due to additional complexities associated with determination, computation and robustness of optimal designs that are not present in the two-factor case. This paper, therefore, is not a mere generalization of our earlier work. It should be noted that for the general case of $2^k$ experiments with binary response, Dorta-Guerra, Gonz\'{a}lez-D\'{a}vila and Ginebra (2008) have obtained an expression for the D-criterion and studied several special cases.

A motivating example is the odor removal study conducted by textile engineers at the University of Georgia. The scientists study the manufacture of bio-plastics from algae that contain odorous volatiles. These odorous volatiles, generated from algae bio-plastics, either occur naturally within the algae or are generated through the thermoplastic processing due to heat and pressure. In order to commercialize these algae bio-plastics, the odor causing volatiles must be removed. Static headspace microextraction and gas chromatography $-$ mass spectroscopy are used to identify the odorous compounds and qualitatively assess whether or not the volatiles have been successfully removed. The outcome of this assessment is the response of the experiment. For that purpose, a study was conducted with a $2^{4-1}_{IV}$ design, a regular fraction, with five replicates using algae and synthetic plastic resin blends. The four different factors were: type of algae, scavenger material (adsorbent), synthetic resin and compatabilizers (see Table~\ref{tab:odorfac} for details).

\begin{table}[ht]\caption{Factors and levels, odor experiment}\label{tab:odorfac}
{\small
\begin{center}
\begin{tabular}{ccp{1.3in}p{1.3in}}
\hline
Factor & Levels & \multicolumn{1}{c}{$-$} & \multicolumn{1}{c}{$+$} \\
\hline
A & algae & raffinated or solvent extracted algae & catfish pond algae \\
B & scavenger material & Aqua Tech activated carbon & BYK-P 4200 purchased from BYK Additives Instruments \\
C & synthetic resin & polyethylene & polypropylene\\
D & compatabilizers & absent & present\\
\hline
\end{tabular}
\end{center}
}
\end{table}

We obtain theoretical results and algorithms for locally optimal designs for $k$ qualitative factors at two levels each and a binary response in the generalized linear model setup. We consider $D$-optimal designs, which maximize the determinant of the information matrix. Although we explore designs for full factorials, i.e., ones in which observations are taken at every possible level combination, when the number of factors is large, full factorials are practically infeasible. Hence the study of fractional factorial designs occupies a substantial part of the linear-model based design literature, and we too study these designs in our setup. A natural question that arises when we use local optimality is whether the resulting designs are robust to the assumed parameter values. We consider this in Section~5.

An alternative approach to design optimality is Bayes optimality (Chaloner and Verdinelli, 1995). For our problem, however, for large $k$ ($k\ge 4$) the computations quickly become expensive. Hence as a surrogate criterion, we explore a D-optimality criterion with the information matrix replaced by its expectation under the prior. This is one of the suggested alternatives to formal Bayes optimality in Atkinson, Donev and Tobias (2007). It has been used by Zayats and Steinberg (2010) for optimal designs for detection capability of networks. We call this {\it EW D-optimality} ($E$ for expectation, $W$ for the notation $w_i$ used for the GLM ``weight'', which can be thought of as information contained in an individual observation). Effectively this reduces to a locally optimal design with local values of the weight parameters replaced by their expectations. The EW D-optimal designs are very good and easy-to-compute surrogates for Bayes D-optimal designs. Unless $k$ is small or the experimenter is quite certain about the parameter values, we suggest the use of EW D-optimal designs. Note that the use of surrogates of Bayes optimality has been recommended by Gotwalt et al.~(2009).

Beyond theoretical results, the question that may be asked is whether these results give the user any advantage in real experiments. It turns out that when $k>2$, in most situations, we gain considerably by taking advantage of the results of this paper instead of using standard linear-model results. Unlike the linear model case, not all nonsingular regular fractions have the same D-efficiency. Indeed, if we have some knowledge of the parameters, we will be able to identify an efficient fractional factorial design, which is often not a regular fraction.

This paper is organized as follows. In Section~2 we describe the preliminary setup. In Section~3 we provide several results for locally D-optimal designs, including the uniqueness of the D-optimal designs, characterization for a design to be locally D-optimal, the concept of EW D-optimal designs, and algorithms for finding D-optimal designs. In Section~4 we discuss the properties of fractional factorial designs. We address the robustness of D-optimal designs in Section~5 and revisit the odor {example in Section~6. Some concluding remarks and topics for future research are discussed in Section~7. Additional results, proofs and some details on the algorithms are relegated to the Supplementary Materials.

%\clearpage

\fontsize{10.95}{14pt plus.8pt minus .6pt}\selectfont
\setcounter{chapter}{2}
\setcounter{equation}{0} %-1

{\bf 2. Preliminary Setup}

Consider a $2^{k}$ experiment with binary response, i.e., an experiment with $k$ explanatory variables at $2$ levels each. Suppose $n_{i}$ units are allocated to the $ i$th experimental condition such that $n_{i}\geqslant 0,$ $i=1,\ldots,2^{k}$, and $n_{1}+\cdots+n_{2^{k}}=n$. We suppose that $n$ is fixed and the problem is to determine the ``optimal" $n_{i}$'s. \ In fact, we write our optimality criterion in terms of the proportions \vspace{-.1in}%
\begin{equation*}
p_{i}=n_{i}/n,\quad i=1,\ldots,2^{k}
\end{equation*}
and determine the ``optimal" $p_i \ge 0$ satisfying $\sum_{i=1}^{2^k}p_i=1$. Since $n_{i}$'s are integers, an optimal design obtained in this fashion may not always be  {viable}. In Section~\ref{algo2} we will consider the design problem over integer $n_i$'s.

We will use a generalized linear model setup. Suppose $\eta$ is a linear predictor that involves the main effects and interactions that are assumed to be in the model. For instance, for a $ 2^{3} $ experiment with a model that includes the main effects and the two-factor interaction of factors $1$ and $2$, $\eta =\beta _{0}+\beta _{1}x_{1}+\beta _{2}x_{2}+\beta _{3}x_{3}+\beta _{12}x_{1}x_{2},$ where each $x_{i}\in \left\{ -1,1\right\}$. The aim of the experiment is to obtain inferences about the parameter vector of factor effects $\boldsymbol\beta =\left( \beta _{0},\beta _{1},\beta _{2},\beta _{3},\beta _{12}\right) ^{\prime }$. In the framework of generalized linear models, the expectation of the response $Y$, $E\left( Y\right) =\pi$, is connected to the linear predictor $\eta$ by the link function $g$:  $\eta =g\left( \pi \right)$ (McCullagh and Nelder, 1989). For a binary response, the commonly used link functions are logit, probit, log-log, and complementary log-log links.

The maximum likelihood estimator of $\boldsymbol\beta $ has an asymptotic covariance matrix (McCullagh and Nelder, 1989; Khuri, Mukherjee, Sinha and Ghosh, 2006) that is the inverse of $ nX^{\prime }WX$, where $W={\rm diag}\left\{ w_{1}p_{1}, ..., w_{2^{k}}p_{2^{k}}\right\} ,$ $w_{i} = \left( \frac{d\pi _{i}}{ d\eta _{i}}\right) ^{2}/( \pi _{i}(1-\pi _{i}) ) \geq 0$, $\eta_i$ and $\pi_i$ correspond to the $i$th experimental condition for $\eta$ and $\pi$, and $X$ is the ``{model matrix}''. For example, for a $2^3$ experiment with model $\eta = \beta_0 + \beta_1 x_1 + \beta_2 x_2 + \beta_3 x_3 + \beta_{12}x_1x_2$,\
\begin{eqnarray}\label{xmat}
X &=& \left(
\begin{array}{ccccc}
+1 & +1 & +1 & +1 & +1\\
+1 & +1 & +1 & -1 & +1\\
+1 & +1 & -1 & +1 & -1\\
+1 & +1 & -1 & -1 & -1\\
+1 & -1 & +1 & +1 & -1\\
+1 & -1 & +1 & -1 & -1\\
+1 & -1 & -1 & +1 & +1\\
+1 & -1 & -1 & -1 & +1\\
\end{array}
\right)
\end{eqnarray}

The $n_i$'s determine how many observations are made at each experimental condition, which are characterized by the rows of $X$. A D-optimal design maximizing $|X'WX|$ depends on the $w_i$'s, which in turn depend on the regression parameters $\boldsymbol\beta$ and the link function $g$. In this paper, we discuss D-optimal designs in terms of $w_i$'s so that our results are not limited to specific link functions.

Unlike experiments with continuous factors, the $2^k$ design points in our setup are fixed and we only have the option of determining the optimal proportions. For results on optimal designs with continuous factors in the GLM setup, see for example, Stufken and Yang (2012).

\fontsize{10.95}{14pt plus.8pt minus .6pt}\selectfont
\setcounter{chapter}{3}
\setcounter{equation}{0} %-1

{\bf 3. Locally D-Optimal Designs}

In this section, we start with a formulation of the local D-optimality problem and establish some general results. Consider a $2^k$ experiment. The goal is to find an optimal ${\mathbf p}=(p_1$, $p_2$, $\ldots$, $p_{2^k})'$ which maximizes $f({\mathbf p}) := |X'WX|$ for specified values of $w_i \geq 0, i=1, \ldots, 2^k$. The specification of the $w_i$'s come from the initial values of the parameters and the link function. Here $p_i\geq 0$, $i=1,\ldots, 2^k$ and $\sum_{i=1}^{2^k} p_i = 1$. It is easy to see that there always exists a D-optimal allocation ${\mathbf p}$ since the set of all feasible allocations is bounded and closed. On the other hand, the uniqueness of D-optimal designs is usually not guaranteed (see Remark~\ref{remarkunique2}). Note that even if all the $p_i$'s are positive, the resulting design is not full factorial in the traditional sense where equal number of replicates are used. On the other hand, if some of the $p_i$'s are zero, then it becomes a fractional factorial design which will be discussed in the next section. In fact, the number of nonzero $p_i$'s in the optimal design could be much less than $2^k$, as we will see in Section~3.3.

\setcounter{section}{1}

{\bf 3.1 Characterization of locally D-optimal designs}

Suppose the parameters (main effects and interactions) are $\boldsymbol\beta = (\beta_0, \beta_1, \ldots,$ $ \beta_d)'$, where $d\geq k$. The following lemma expresses the objective function as an order-($d$+1) homogeneous polynomial of $p_1, \ldots, p_{2^k}$~.

\begin{lemma}\label{xwx}
Let  $X[i_1,i_2,\ldots,i_{d+1}]$ be the $(d+1)\times (d+1)$ sub-matrix consisting of the $i_1\mbox{th}$, $i_2\mbox{th}$, $\ldots$, $i_{d+1}\mbox{th}$ rows of the {model matrix} $X$. Then
$$f({\mathbf p})=|X'WX|=\sum_{1\leq i_1<\cdots<i_{d+1}\leq 2^k} |X[i_1,i_2,\ldots,i_{d+1}]|^2 \cdot p_{i_1}w_{i_1}p_{i_2}w_{i_2}\cdots p_{i_{d+1}}w_{i_{d+1}} .$$
\end{lemma}

Gonz\'{a}lez-D\'{a}vila, Dorta-Guerra and Ginebra (2007, Proposition~2.1) obtained essentially the same result. This can also be proved directly using the results from Rao (1973, Chapter 1). From Lemma~\ref{xwx} it is immediate that at least $(d+1)$ $w_i$'s, as well as the corresponding $p_i$'s, have to be positive for the determinant $f({\mathbf p})$ to be nonzero. This implies that if ${\mathbf p}$ is D-optimal, then $p_i < 1$ for each $i$. Theorem~\ref{theorem30} below gives a sharper bound, $p_i \leq \frac{1}{d+1}$ for each $i=1, \ldots, 2^k$, for the optimal allocation. Let us define for each $i=1, \ldots, 2^k$,
\begin{equation}\label{f_i(x)}
f_i(z) =f\left(\frac{1-z}{1-p_i}p_1,\ldots,\frac{1-z}{1-p_i}p_{i-1},z, \frac{1-z}{1-p_i}p_{i+1},\ldots, \frac{1-z}{1-p_i}p_{2^k}\right), \>\>\> 0\leq z\leq 1.
\end{equation}
Note that $f_i(z)$ is well defined for all ${\mathbf p}$ of interest (that is, $p_i < 1$ for each $i$).

\begin{theorem}\label{theorem30}
Suppose $f\left({\mathbf p}\right)>0$. Then ${\mathbf p}$ is D-optimal if and only if for each $i=1, \ldots, 2^k$, one of the two conditions below is satisfied:
\begin{itemize}
\item[(i)] $p_i=0$ and $f_i\left(\frac{1}{2}\right) \leq\frac{d+2}{2^{d+1}}f({\mathbf p})$;
\item[(ii)] $0 < p_i \leq \frac{1}{d+1}$ and $f_i(0) =\frac{1-p_i(d+1)}{(1-p_i)^{d+1}}f({\mathbf p})$.
\end{itemize}
\end{theorem}

\begin{remark}\label{remarkunique1}
{\rm
Theorem~\ref{theorem30} is essentially a specialized version of the general equivalence theorem on a pre-determined finite set of design points. Unlike the usual form of the equivalence conditions (for examples, see Kiefer (1974), Pukelsheim (1993), Atkinson et al.~(2007), Stufken and Yang (2012), Fedorov and Leonov (2014)) where the inverse matrix of $X'WX$ needs to be calculated, Theorem~\ref{theorem30} is expressed in terms of the determinant quantities $f({\mathbf p})$, $f_i(\frac{1}{2})$ and $f_i(0)$ only. These expressions are critical for the algorithms proposed later in this section. This theorem also gives a sharper bound $0<p_i\leq 1/(d+1)$ for support points. Note that even if $p_i=0$ for some $i$, it is still possible that the equality $f_i(1/2)=(d+2)/(2^{d+1})\cdot f({\mathbf p})$ holds. In the {Supplementary Materials}, we provide a self-contained proof of Theorem~\ref{theorem30} which does not rely on any general equivalence theorem. Its connection to the General Equivalence Theorem is provided in the Supplementary Materials.
}\end{remark}

Designs that are supported on $(d+1)$ points are attractive in many experiments because they require a minimum number of settings. In our context, a design ${\mathbf p} = (p_1, \ldots, p_{2^k})'$ is called {\it minimally supported} if it has exactly $(d+1)$ nonzero $p_i$'s. For designs supported on rows $i_1,\ldots, i_{d+1}$, the $D$-optimal choice of weights is $p_{i_1} = \cdots = p_{i_{d+1}} = 1/(d + 1)$. This result can be obtained from Lemma~\ref{xwx} directly. Yang et al. (2012) found a necessary and sufficient condition for a minimally supported design to be D-optimal for $2^2$ main-effects model. With the aid of Theorem~\ref{theorem30}, we provide a generalization for $2^k$ designs in the next theorem. Note that $w_i > 0$ for each $i$ for the  commonly used link functions including logit, probit, and (complementary) log-log.

\begin{theorem}\label{dsaturatedconditiontheorem}
Assume $w_i > 0$, $i=1, \ldots, 2^k$. Let ${\mathbf I} = \{i_1,\ldots,i_{d+1}\}\subset\{1,\ldots,2^k\}$ be an index set satisfying $|X[i_1,\ldots,$ $i_{d+1}]|\neq 0$. Then the minimally supported design satisfying $p_{i_1}=p_{i_2}=\cdots=p_{i_{d+1}}=\frac{1}{d+1}$ is D-optimal if and only if for each $i\notin {\mathbf I}$,
$$
\sum_{j\in {\mathbf I}} \frac{|X[\{i\}\cup {\mathbf I}\setminus\{j\}]|^2}{w_j} \leq \frac{|X[i_1,i_2,\ldots,i_{d+1}]|^2}{w_i}.
$$
\end{theorem}

For example, under the $2^2$ main-effects model, since $|X[i_1, i_2, i_3]|^2$ is constant across all choices of $i_1,i_2,i_3$, $p_1=p_2=p_3=1/3$ is D-optimal if and only if $v_1+v_2+v_3\leq v_4$, where $v_i = 1/w_i$, $i=1,2,3,4$. This gives us Theorem~1 of Yang, Mandal and Majumdar (2012). For the $2^3$ main-effects model, the {model matrix} $X$ is given by (\ref{xmat}) with the last column deleted. Using this order of rows, the standard regular fractional factorial design $p_1=p_4=p_6=p_7=1/4$ given by the defining relation $1=ABC$ is D-optimal if and only if $v_1+v_4+v_6+v_7 \leq 4\min\{v_2,v_3,v_5,v_8\}$, and the other standard regular fractional design $p_2=p_3=p_5=p_8=1/4$ is D-optimal if and only if $v_2+v_3+v_5+v_8 \leq 4\min\{v_1,v_4,v_6,v_7\}$.

\begin{remark}\label{remarkunique2}
{\rm
In order to characterize the uniqueness of the optimal allocation, we define a matrix $X_w =[{\mathbf 1},\ {\mathbf w}*{\mathbf 1},\ {\mathbf w}*\gamma_2,\ \ldots,\ {\mathbf w}*\gamma_s]$, where ${\mathbf 1}$ is the $2^k\times 1$ vector of all $1$'s, $\{{\mathbf 1}, \gamma_2, \ldots, \gamma_s\}$ forms the set of all distinct pairwise Schur products {(or entrywise product)} of the columns of the {model matrix} $X$, ${\mathbf w} = (w_1, \ldots, w_{2^k})'$, and ``$*$" indicates Schur product. It can be verified that any two feasible allocations ($p_i \ge 0$ satisfying $\sum_{i=1}^{2^k}p_i = 1$) generate the same matrix $X'WX$ as long as the difference of the matrices belongs to the null space of $X_w$. If rank($X_w$) $< 2^k$, any criterion based on $X'WX$ yields an affine set of solutions with dimension $2^k-{\rm rank}(X_w )$. If rank($X_w $) $= 2^k$, the D-optimal allocation ${\mathbf p}$ is unique. For example, for a $2^3$ design the model consisting of all main effects and one two-factor interaction, or for a $2^4$ design the model consisting of all main effects, all two-factor interactions, and one three-factor interaction, the D-optimal allocation is unique.
}\end{remark}

\setcounter{section}{2}

  {\bf 3.2 EW D-optimal designs}\label{sectionew}

Since locally D-optimal designs depend on $w_i$'s, they require assumed values of $w_i$'s, or $\beta_i$'s, as input. In Section~5, we will examine the robustness of D-optimal designs to mis-specification of $\beta_i$'s. An alternate to local optimality is Bayes optimality (Chaloner and Verdinelli, 1995). In our setup, a Bayes D-optimal design maximizes $E(\log|X'WX|)$ where the expectation is taken over the prior on $\beta_i$'s. One difficulty of Bayes optimality is that it is computationally expensive. In order to overcome this drawback we explore an alternative suggested by Atkinson, Donev and Tobias (2007) where $W$ in the Bayes criterion is replaced by its expectation. We call this {\it EW} (expectation of $W$) D-optimality.

{\bf Definition:} An {\it EW D-optimal} design is an optimal allocation ${\mathbf p}$ that maximizes $|X'E(W)X|$.

Note that EW D-optimality may be viewed as local D-optimality with $w_i$'s replaced by their expectations. All of the existence and uniqueness properties of locally D-optimal design apply. Since $w_i > 0$ for all $\boldsymbol\beta$ under typical link functions, $E(w_i)>0$ for each $i$. By Jensen's inequality, $$E\left(\log|X'WX|\right) \leq \log|X' E(W) X|$$ since $\log|X^\prime W X|$ is concave in $\bf w$. Thus an EW D-optimal design maximizes an upper bound for Bayesian D-optimality criterion.

In practice, once $E(w_i)$'s are calculated via numerical integration, algorithms for local D-optimality can be applied with $w_i$ replaced by $E(w_i)$. We will show that EW D-optimal designs are often almost as efficient as designs that are optimal with respect to the Bayes D-optimality criterion, while realizing considerable savings in computation time.  {In fact, while searching for a EW D-optimal design, the integration can be performed in advance of the optimization. This provides a computational advantage over the search for Bayesian D-optimal designs, where integration needs to be performed in each step of the optimization, in order to evaluate the design.} Furthermore, EW D-optimal designs are highly robust in terms of maximum loss of efficiency (Section~5).

Given link function $g$, let $\nu =  \left[\left(g^{-1}\right)'\right]^2/\left[g^{-1}(1-g^{-1})\right]$. Then $w_i = \nu(\eta_i) = \nu\left({\mathbf x_i}'\boldsymbol\beta\right)$, $i=1, \ldots, 2^k$, where ${\mathbf x}_i$ is the $i$th row of the {model matrix} $X$, and $\boldsymbol{\beta} = (\beta_0, \beta_1, \ldots,$ $\beta_d)'$. Suppose the regression coefficients $\beta_0,\beta_1,\ldots,\beta_d$ are independent, and $\beta_1, \ldots, \beta_d$ each has a symmetric distribution about $0$ (not necessarily the same distribution), then all the $w_i,\ i=1, \ldots, 2^k$ have the same distribution and the uniform design $p_1 = \cdots = p_{2^k} = 2^{-k}$ is an EW D-optimal design for any given link function (by ``uniform design'' we mean a design with uniform allocation on its support points). On the other hand, in many experiments we may be able to assume that the slope of a main effect is non-decreasing. If $\beta_i \in [0, \beta_{iu}]$ for each $i$, the uniform design will not be EW D-optimal in general, as illustrated in the following example.

\begin{example}\label{23mainnonuniform}{\rm
Consider a $2^3$ experiment with main-effects model. Suppose $\beta_0$, $\beta_1$, $\beta_2$ and $\beta_3$ are independent, $\beta_0 \sim U[-3, 3]$, and $\beta_1, \beta_2, \beta_3 \sim U[0, 3]$. Then $E(w_1)=E(w_8)=0.042$, $E(w_2) = E(w_3)= \cdots = E(w_7) = 0.119$. Under the logit link the EW D-optimal design is ${\mathbf p}_{e} = (0, 1/6, 1/6, 1/6,$ $1/6,$ $1/6,$ $1/6,$ $0)'$, and the Bayesian D-optimal design, which maximizes $\phi({\mathbf p})$ $=$ $E(\log|X'WX|)$, is ${\mathbf p}_o$ = ($0.004$, $0.165$, $0.166$, $0.165$, $0.165$, $0.166$, $0.165$, $0.004$)$'$. The efficiency of ${\mathbf p}_{e}$ with respect to ${\mathbf p}_o$ is $\exp\left\{\frac{\phi({\mathbf p}_e)-\phi({\mathbf p}_o)}{d+1}\right\} \times 100\% = 99.98\%$, while the efficiency of the uniform design is $94.39\%$. Also note, in this example, the EW and Bayes criteria lead to virtually the same design. It is remarkable that it takes $2.39$ seconds to find an EW solution while it takes $121.73$ seconds to find a Bayes solution. The difference in computational time is even more prominent for $2^4$ case (24 seconds versus 3147 seconds). All multiple integrals here are calculated using  {\tt R} function {\tt adaptIntegrate} in the package {\tt cubature}.
}
\end{example}

\setcounter{section}{3}

  {\bf 3.3 Algorithms to search for locally D-optimal allocation}\label{section32}

  In this section, we develop efficient algorithms to search for locally D-optimal allocations with given $w_i$'s. The same algorithms can be used for finding EW D-optimal designs.

\subsection{Lift-one algorithm for maximizing $f({\mathbf p})=|X'WX|$}\label{algo1}

 Here we propose the {\it lift-one} algorithm for obtaining locally D-optimal ${\mathbf p} = (p_1, \ldots, p_{2^k})'$ with given $w_i$'s. The basic idea is that, for randomly chosen $i \in \{1, \ldots, 2^k\}$, we update $p_i$ to $p_i^*$ and all the other $p_j$'s to $p_j^*=p_j\cdot \ \frac{1-p_i^*}{1-p_i}$. {This technique is motivated by the coordinate descent algorithm (Zangwill, 1969).} It is also in spirit similar to the idea of one-point correction in the literature (Wynn, 1970; Fedorov, 1972; M\"uller, 2007), where design points are added/adjusted one by one. The major advantage of the lift-one algorithm is that in order to determine an optimal $p_i^*$, we need to calculate $|X'WX|$ only once due to Lemma~\ref{xwx} {(see Step~3$^\circ$ of the algorithm below)}.

 {\it Lift-one algorithm:}
\begin{itemize}
 \item[$1^\circ$] Start with arbitrary ${\mathbf p}_0=(p_1,\ldots,p_{2^k})'$ satisfying $0<p_i<1$, $i=1,\ldots,2^k$ and compute $f\left({\mathbf p}_0\right)$.
 \item[$2^\circ$] Set up a random order of $i$ going through $\{1,2,\ldots,2^k\}$.
 \item[$3^\circ$] Following the random order of $i$ in 2$^\circ$, for each $i$, determine $f_i(z)$ as in (\ref{lem:f_i(x)}) in Supplementary Materials. In this step, either $f_i(0)$ or $f_i\left(\frac{1}{2}\right)$ needs to be calculated according to equation~(\ref{f_i(x)}).
 \item[$4^\circ$] Define ${\mathbf p}_*^{(i)} =\left(\frac{1-z_*}{1-p_i}p_1,\ldots,\frac{1-z_*}{1-p_i}p_{i-1},z_*, \frac{1-z_*}{1-p_i}p_{i+1},\ldots, \frac{1-z_*}{1-p_i}p_{2^k}\right)'$, where $z_*$ maximizes $f_i(z)$ with $0\leq z\leq 1$ (see Lemma~\ref{algo1lemma31}). Note that $f({\mathbf p}_*^{(i)})=f_i(z_*)$. Lemma~\ref{algo1lemma31} gives a simple analytical formula for the update in terms of $f_i(0)$ or $f_i(1/2)$.
 \item[$5^\circ$] Replace ${\mathbf p}_0$ with ${\mathbf p}_*^{(i)}$, $f\left({\mathbf p}_0\right)$ with $f({\mathbf p}_*^{(i)})$.
 \item[$6^\circ$] Repeat $2^\circ\sim 5^\circ$ until convergence, that is, $f({\mathbf p}_0)=f({\mathbf p}_*^{(i)})$ for each $i$.
\end{itemize}

While in all examples that we studied, the lift-one algorithm converges very fast, we do not have a proof of convergence. There is a modified lift-one algorithm, which is only slightly slower, that can be shown to converge. This algorithm can be described as follows. For the $10m$th iteration and a fixed order of $i=1, \ldots, 2^k$ we repeat steps $3^\circ\sim 5^\circ$, $m=1,2, \ldots$. If ${\mathbf p}^{(i)}_*$ is a better allocation found by the lift-one algorithm than the allocation ${\mathbf p}_0$, instead of updating ${\mathbf p}_0$ to ${\mathbf p}^{(i)}_*$ immediately, we obtain ${\mathbf p}^{(i)}_*$ for each $i$, and replace ${\mathbf p}_0$ with the first best one among $\left\{ {\mathbf p}^{(i)}_*, i=1, \ldots, 2^k\right\}$. It should be noted that the updating strategy at the $10m$th iteration here is similar to the Fedorov-Wynn algorithm (Fedorov (1972), Fedorov and Hackl (1997)) but with a more efficient updating formula. For iterations other than the $10m$th, we follow the original lift-one algorithm update.

\begin{theorem}\label{algo1theorem15}
When the lift-one algorithm or the modified lift-one algorithm converges, the resulting allocation ${\mathbf p}$ maximizes $|X'WX|$ on the set of feasible allocations. Furthermore, the modified lift-one algorithm is guaranteed to converge.
\end{theorem}

%\bigskip
\begin{figure}[h]
\centering
\makebox{\includegraphics[scale=.23]{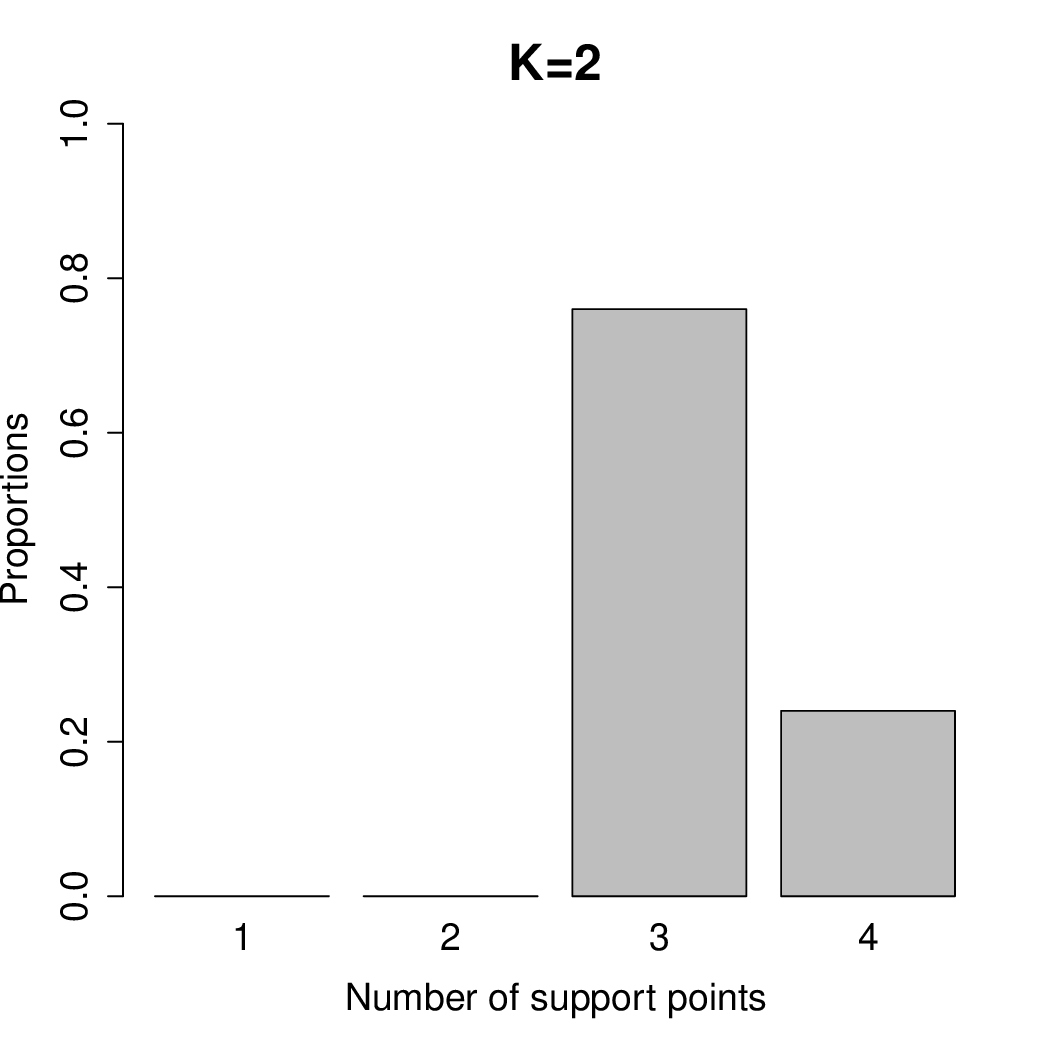}}
\makebox{\includegraphics[scale=.23]{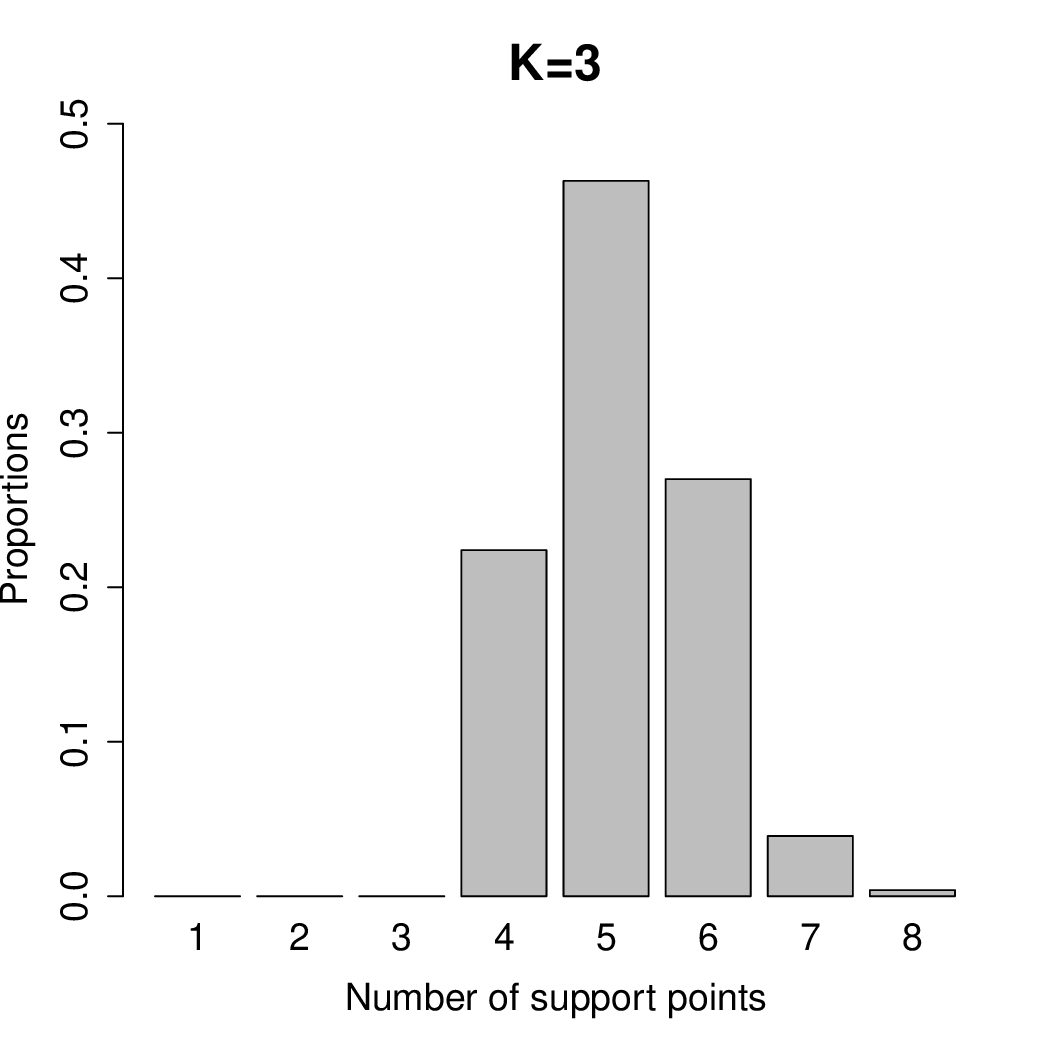}}
\makebox{\includegraphics[scale=.23]{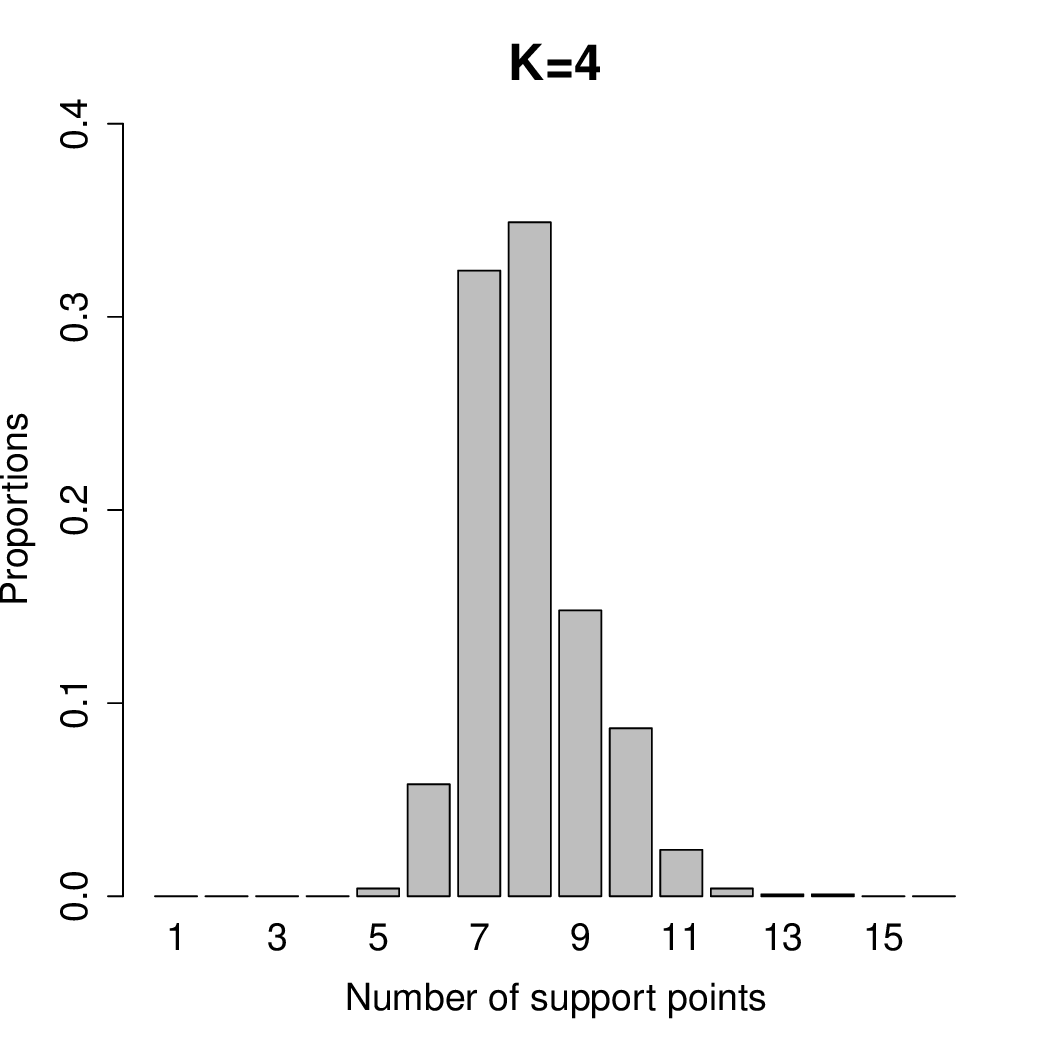}}
\makebox{\includegraphics[scale=.23]{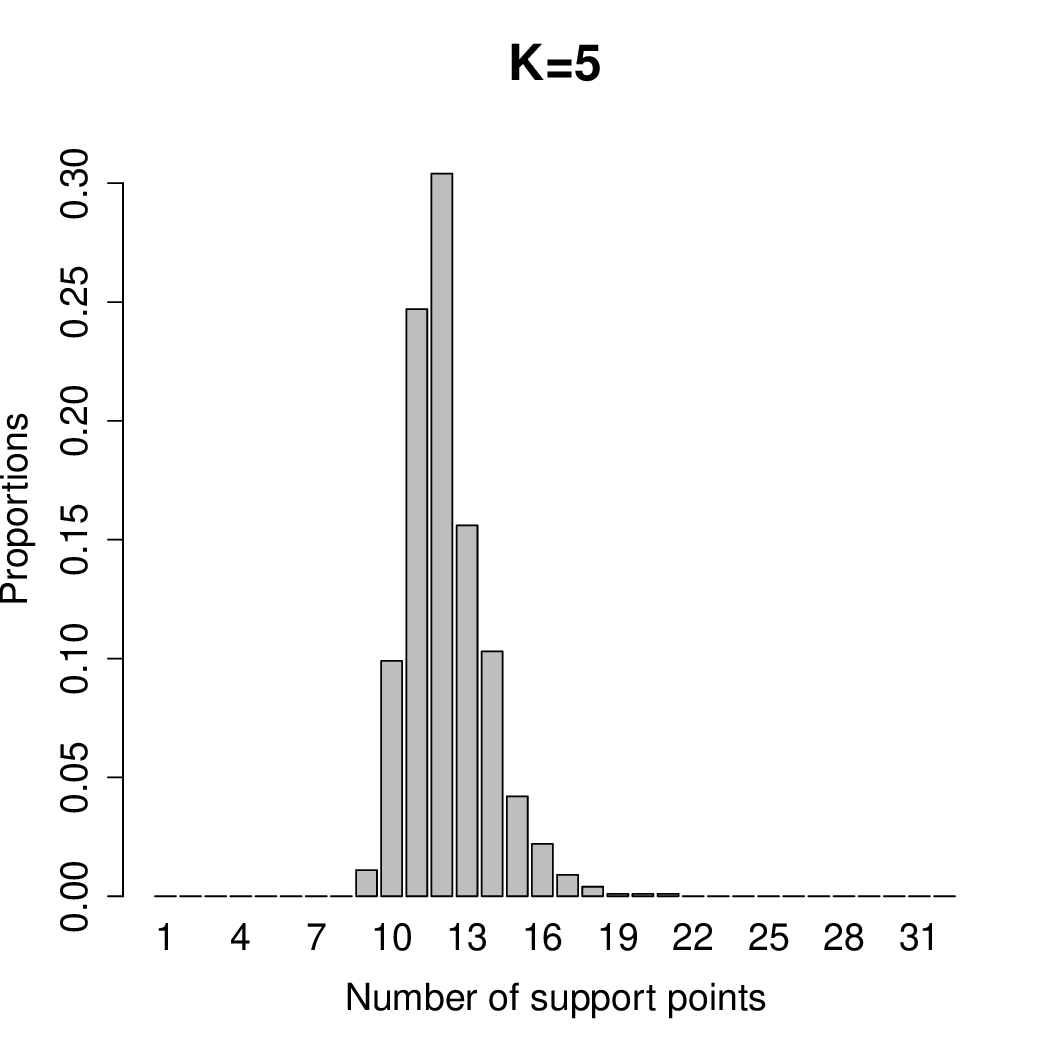}}
\makebox{\includegraphics[scale=.23]{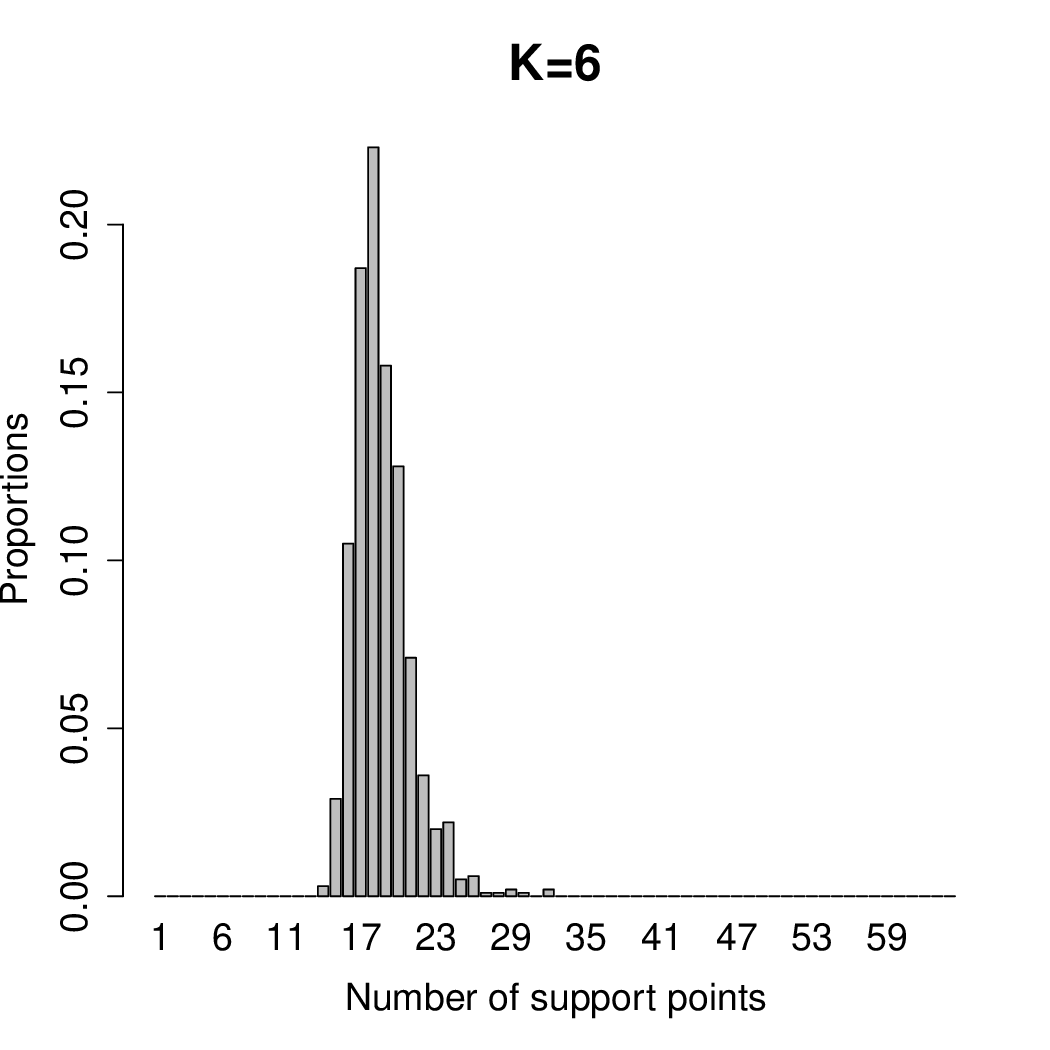}}
\makebox{\includegraphics[scale=.23]{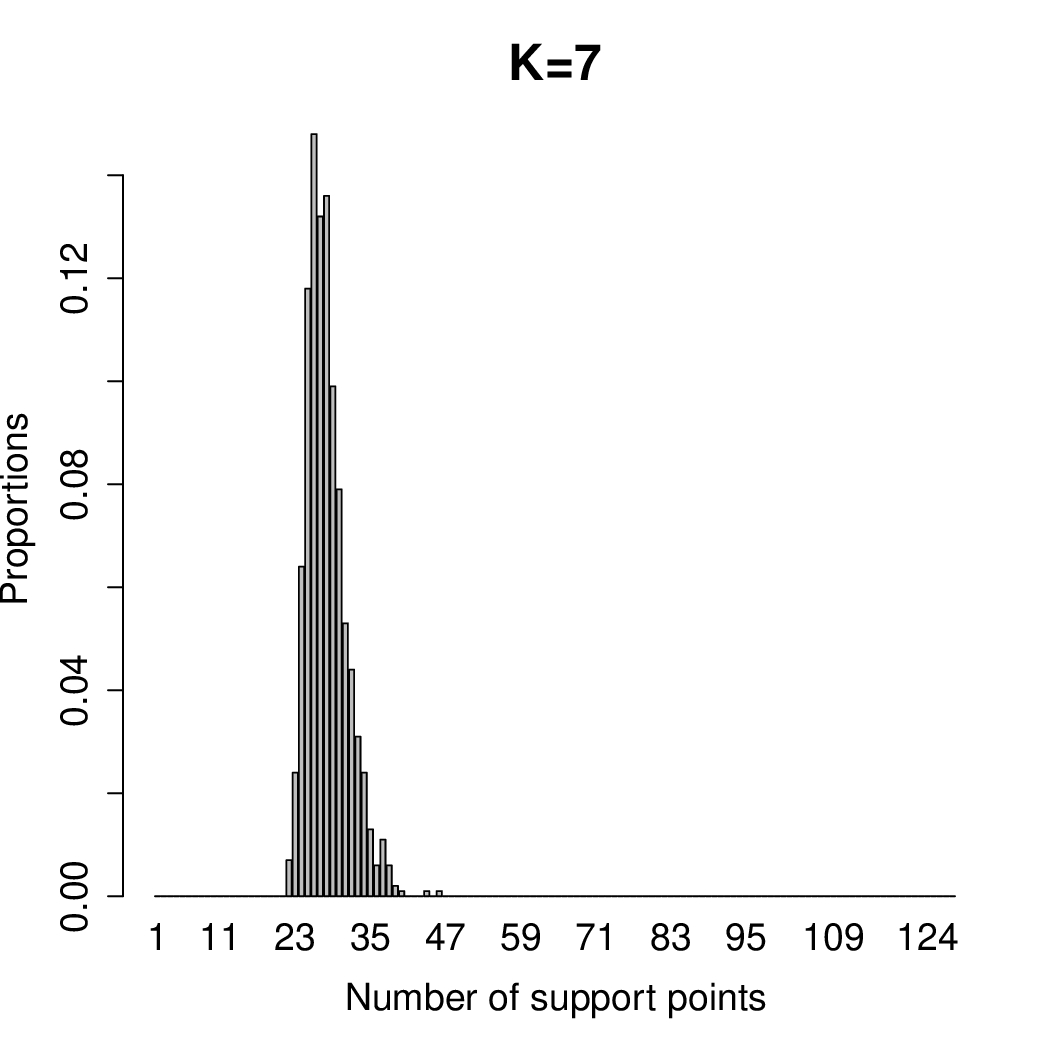}}
%\vskip -0.5in
\caption{\label{figsim}Number of support points in an optimal design (based on 1000 simulations)}
\end{figure}

Our simulation studies indicate that as $k$ grows, the optimal designs produced by the lift-one algorithm for main-effects models is supported only on a fraction of all the $2^k$ design points. To illustrate this, we randomly generate the regression coefficients i.i.d.~from $U(-3,3)$ and apply our algorithm to find the optimal designs under the logit link. Figure~\ref{figsim} gives histograms of numbers of support points in optimal designs found by the lift-one algorithm. For example, with $k=2$, 76\% of the designs are supported on three points only and 24\% of them are supported on all four points. As $k$ becomes larger, the number of support points moves towards a smaller fraction of $2^k$. On the other hand, a narrower range of coefficients requires a larger portion of support points. For example, the mean numbers of support points with $\beta_i$'s i.i.d.~from $U(-3,3)$ are $3.2, 5.1, 8.0, 12.4, 18.7, 28.2$ for $k=2, 3, 4, 5, 6, 7$, respectively. The corresponding numbers increase to $4.0, 7.1, 11.9, 19.1, 30.6, 47.7$ for $U(-1,1)$, and further to $4.0, 7.6, 14.1, 24.7, 41.2, 66.8$ for $U(-0.5, 0.5)$.

The lift-one algorithm is much faster than commonly used optimization techniques (Table~\ref{tab:cpu}) including Nelder-Mead, quasi-Newton, conjugate-gradient, simulated annealing  (for a comprehensive reference, see Nocedal and Wright (1999)), as well as popular design algorithms for similar purposes including Fedorov-Wynn (Fedorov (1972), Fedorov and Hackl (1997),  Fedorov and Leonov (2014)), Multiplicative (Titterington (1976, 1978), Silvey et al.~(1978)), and Cocktail (Yu (2010)) algorithms. We utilize the function {\tt constrOptim} in {\tt R} to implement Nelder-Mead, quasi-Newton, conjugate-gradient, and simulated annealing algorithms. As the number of design points ($2^k$) increases, those four algorithms fail to achieve adequately accurate solutions (marked by ``$-$" in Table~\ref{tab:cpu} which indicates that the relative efficiency compared with the lift-one solutions is below 80\% on average). For example, it takes the Nelder-Mead algorithm 51.73 seconds to find solutions (${\mathbf p}_{NM}$) at $k=6$ whose relative efficiency (defined as $\left(f({\mathbf p}_{NM})/f({\mathbf p}_{lo})\right)^{1/(k+1)}$) compared with the lift-one solutions (${\mathbf p}_{lo}$) is only 65\% on average. As $k$ increases from $2$ to $3$, although the time spent for simulated annealing algorithm reduces from 83.09 seconds to 18.54 seconds, the relative efficiency on average decreases from 99.8\% to 93.0\% (it drops down to 66\% at k=4 and 54\% at k=5). The relative efficiencies do not improve much if more iterations or multiple initial points are allowed. The implementation of the Fedorov-Wynn algorithm here is mainly based on Fedorov and Leonov (2014, \S 3.1) with updating formula for $(X'WX)^{-1}$. As for the Multiplicative and Cocktail algorithms, we followed Yu (2010) and Mandal, Wong and Yu (2014). Each of these three algorithms achieves essentially the same efficiency compared to the lift-one algorithm. For a fair comparison, all the programs were written in {\tt R}, controlled by the same relative convergence tolerance $10^{-5}$, and run at the same computer with Intel CPU at 2.5GHz, 8GB memory, and 64-bit (Windows 8.1) Operating System. Based on the simulation results shown in Table~\ref{tab:cpu}, the lift-one algorithm runs at a much faster speed across different model setups. In terms of the number of support points on average, only the solutions found by the Cocktail algorithm are comparable with lift-one solutions. Typically, the Multiplicative algorithm finds twice as many support points as lift-one's, while the other five algorithms simply keep positive weights on all the $2^k$ design points.

\begin{table}[h]\caption{\label{tab:cpu}Performance of the lift-one algorithm (CPU time in seconds for 100 simulated {${\boldsymbol\beta}$ from U($-3,3$) with logit link and main-effects model})}
\centering
{\footnotesize
\fbox{%
{
\begin{tabular}{c|rrrrrrrr}
  & \multicolumn{8}{c}{Algorithms}\\
 Designs   & Nelder-      & quasi-   & conjugate& simulated & Fedorov  & Multipli- & Cocktail &  Proposed \\
           & Mead         & Newton   & gradient & annealing & -Wynn    & cative    &          & lift-one  \\
\hline
$2^2$      &       1.42   &    0.19   &     2.09   & 83.09    &    6.14  &     0.28  &    0.16  &  0.11    \\
$2^3$      &      8.76   &   24.64  &   171.74  &  18.54    &   11.25  &     0.86  &    0.53  &  0.36    \\
$2^4$      &     17.88   &   $-$       &    $-$       &  $-$        &   21.77  &    10.97  &    4.46  &  1.07    \\
$2^5$      &     31.64         &    $-$      &    $-$       &   $-$       &   47.66  &    50.12  &   68.88  &  4.82    \\
$2^6$      &      $-$       &   $-$       &   $-$        &  $-$        &  106.89  &   229.17  &  189.83  & 18.29    \\
$2^7$      &      $-$        &   $-$       &   $-$        &  $-$        &  241.80  &   890.44  &  439.55  & 75.58    %\\[3pt]
 \end{tabular}}
}}
\end{table}

\begin{remark}\label{remarkunique3}
{\rm
{
There are at least two advantages of the proposed algorithm over the competitors listed above. Firstly, the lift-one algorithm exploits the convex structure of the optimization problem (the set of design measures over $\{-1,1\}^k$ is convex, and the objective function $f({\mathbf p})$ is log-concave), whereas  {some of} the other algorithms compared do not. Secondly, Lemma~\ref{xwx} has been used to reduce the number of determinant calculations required per iteration of the algorithm. In Table~\ref{tab:cpu} the comparison with a Federov-Wynn algorithm demonstrates that the gain in speed due to these features of the new algorithm is significant.
}
}\end{remark}

\subsection{Algorithm for maximizing $|X'WX|$ with integer solutions}\label{algo2}

To maximize $|X'WX|$, an alternative algorithm, called {\it exchange algorithm}, is to adjust $p_i$ and $p_j$ simultaneously for randomly chosen index pair $(i,j)$ (see Supplementary Materials for detailed description). The original idea of exchange was suggested by Fedorov (1972). It follows from Lemma~\ref{xwx} that the optimal adjusted $(p_i^*, p_j^*)$ can be obtained easily by maximizing a quadratic function. Unlike the lift-one algorithm, the exchange algorithm can be applied to search for integer-valued optimal allocation ${\mathbf n} = (n_1, \ldots, n_{2^k})'$, where $\sum_i n_i = n$.

 {\it {Exchange algorithm for integer-valued allocations:}}
\begin{itemize}
 \item[$1^\circ$] Start with initial design ${\mathbf n}=(n_1,\ldots,n_{2^k})'$ such that $f({\mathbf n})>0$.
 \item[$2^\circ$] Set up a random order of $(i,j)$ going through all pairs $$\{(1,2),(1,3),\ldots,(1,2^k), (2,3), \ldots,(2^k-1,2^k)\}$$
 \item[$3^\circ$] For each $(i,j)$, let $m=n_i+n_j$. If $m=0$, let ${\mathbf n}^*_{ij}={\mathbf n}$. Otherwise, calculate $f_{ij}(z)$ as given in equation~(\ref{lem:67fij}). Then let
 $${\mathbf n}^*_{ij} =\left(n_1,\ldots,n_{i-1},z_*,n_{i+1},\ldots,n_{j-1},m-z_*,n_{j+1},\ldots,n_{2^k}\right) $$
 where the integer $z_*$ maximizes $f_{ij}(z)$ with $0\leq z\leq m$ according to  Lemma~\ref{algo3lemma1} in Supplementary Materials. Note that $f({\mathbf n}^*_{ij})=f_{ij}(z_*)\geq f({\mathbf n})>0$.
 \item[$4^\circ$] Repeat $2^\circ\sim 3^\circ$ until convergence (no more increase in terms of $f({\mathbf n})$ by any pairwise adjustment).
\end{itemize}

As expected, the integer-valued optimal allocation $(n_1, \ldots, n_{2^k})'$ is consistent with the proportion-valued allocation $(p_1, \ldots, p_{2^k})'$ for large $n$. For small $n$, the algorithm may be used for the fractional design problem in Section~4. It should be noted that the exchange algorithm for integer-valued solutions is not guaranteed to converge to the optimal solutions, especially when $n$ is small compared to $2^k$. However, when we search for optimal proportions, our algorithm with slight modification is guaranteed to converge (see Supplementary Materials for details).

In terms of finding optimal proportions, the exchange algorithm produces essentially the same results as the lift-one algorithm, although the former is relatively slower. For example, based on 1000 {simulated ${\boldsymbol\beta}$'s from U(-3,3) with logit link and the main-effects model}, the ratio of computational time of the exchange algorithm over the lift-one algorithm is 6.2, 10.2, 16.8, 28.8, 39.5 and 51.3 for $k=2,\ldots,7$ respectively. Note that it requires 2.02, 5.38, 19.2, 84.3, 352, and 1245 seconds respectively to finish the 1000 simulations using the lift-one algorithm on a regular PC with 2.26GHz CPU and 2.0G memory. As the total number of factors $k$ becomes large, the computation is more intensive.

 {It should be noted that the general purpose optimization algorithms might be a little slow and faster alternatives should exist. For example, the adaptive barrier method might be inefficient compared to transformations to obtain an unconstrained optimization problem. For the pseudo-Bayesian designs, it is possible that a fixed quadrature scheme would be faster, though possibly less accurate. Detailed study of the computational properties of the proposed algorithms is a topic for future research.}

%\clearpage

\fontsize{10.95}{14pt plus.8pt minus .6pt}\selectfont
\setcounter{chapter}{4}
\setcounter{equation}{0} %-1
\setcounter{section}{1}

  {\bf 4. Fractional Factorial Designs}

If for the optimal allocation some $p_i$'s are zero, then the resulting design is necessarily a fractional factorial one. Even if all of the proportions in the optimal design are substantially away from zero, the experimenter may need, or prefer, to use a fractional factorial design, because even for moderately large values of $k$, the total number of observations $n$ would have to be large to get integer $np_i$'s. For linear models, the accepted practice is to use regular fractions due to the many desirable properties like minimum aberration and optimality. We will show that in our setup the regular fractions are often not optimal. As a first step, however, we start by identifying situations when they are optimal.

We use $2^3$ designs for illustration. The {model matrix} for $2^3$ main-effects model consists of the first four columns of $X$ given in (\ref{xmat}) and $w_j$ represents the information in the $j$th experimental condition, i.e., the $j$th row of $X$. Suppose the maximum number of experimental conditions is fixed at a number less than $8$, and the problem is to identify the experimental conditions and corresponding $p_i$'s that optimize the objective function. Half fractions use $4$ experimental conditions (hence the design is uniform). The half fractions defined by rows $\{1, 4, 6, 7\}$ and $\{2, 3, 5, 8\}$ are regular fractions, given by the defining relations $1=ABC$ and $-1=ABC$ respectively. If all regression coefficients except the intercept are zeros, then the regular fractions are D-optimal, since all the $w_i$'s are equal.  The following Theorem identifies the necessary and sufficient conditions for regular fractions to be D-optimal in terms of $w_i$'s.

\begin{theorem}\label{theorem2^3}
For the $2^3$ main-effects model, suppose $\beta_1=0$ (which implies $w_1=w_5$, $w_2=w_6$, $w_3=w_7$, and $w_4=w_8$). The regular fractions $\{1,4,6,7\}$ and $\{2,3,5,8\}$ are D-optimal {within the class of half-fractions} if and only if
  $$4\ \min\{w_1,w_2,w_3,w_4\}\geq \max\{w_1,w_2,w_3,w_4\}.$$
Suppose $\beta_1=\beta_2=0$ (thus $w_1=w_3=w_5=w_7$ and $w_2=w_4=w_6=w_8$). The two regular half-fractions $\{1,4,6,7\}$ and $\{2,3,5,8\}$ are D-optimal half-fractions if and only if $4 \min\{w_1,w_2\} \geq \max\{w_1,w_2\}.$
\end{theorem}

\begin{example}\label{example2^3}
{\rm
Under logit link, consider the $2^3$ main-effects model with $\beta_1=\beta_2=0$, {which implies $w_1=w_3=w_5=w_7$ and $w_2=w_4=w_6=w_8$}. The regular half-fractions $\{1,4,6,7\}$ and $\{2,3,5,8\}$ {have the same $|X'WX|$ but not the same $X'WX$. They} are D-optimal half-fractions if and only if one of the following happens:
\begin{eqnarray}\label{eq1}
    &(i)& |\beta_3| \leq  \log 2 \\
    &(ii)& |\beta_3| > \log 2 \mbox{\hspace{.1in} and \hspace{.1in} } |\beta_0| \leq \log\left(\frac{2 e^{|\beta_3|} - 1}{e^{|\beta_3|} - 2}\right).\nonumber
\end{eqnarray}
When the regular half-fractions are not optimal, it follows from Lemma~\ref{xwx} that the goal is to find $\{i_1,i_2,i_3,i_4\}$ that maximizes $|X[i_1,i_2,i_3,i_4]|^2w_{i_1}w_{i_2}w_{i_3}w_{i_4}$~. Recall that in this case there are only two distinct $w_i$'s.  If $\beta_0\beta_3 > 0$, $w_i$'s corresponding to $\{2,4,6,8\}$ are larger than others, so this fraction given by $C=-1$ will maximize $w_{i_1}w_{i_2}w_{i_3}w_{i_4}$~. But this leads to a singular {model matrix}. It is not surprising that the D-optimal half-fractions are ``close'' to the design $\{2,4,6,8\}$, and are in fact given by the 16 designs each consisting of three elements from $\{2,4,6,8\}$ and one from $\{1,3,5,7\}$. We call these {\it modified $C=-1$} fractions. All the 16 designs lead to the same $|X'WX|$, which is $w_1w_2^3/4$. For $\beta_0\beta_3 < 0$, D-optimal half-fractions are similarly obtained from the fraction $C=+1$.

Figure~\ref{fig0} partitions the parameter space for $2^3$ main-effects logit model. The left panel corresponds to the case (a) $\beta_1=\beta_2=0$. Here the parameters in the middle region would make the regular fractions D-optimal, whereas the top-right and bottom-left regions correspond to the case $\beta_0\beta_3>0$. Similarly the other two regions correspond to the case $\beta_0\beta_3<0$ so that modified $C=-1$ is optimal. The right panel of Figure~\ref{fig0} is for the case (b) $\beta_1=0$ and shows the contour plots for the largest $|\beta_0|$'s that would make the regular fractions D-optimal. (For details, see Supplementary Materials of this paper.) Along with Figure~\ref{fig0}, conditions~(\ref{eq1}) and (\ref{eq:2^3b1=0}) in Supplementary Materials indicate that if $\beta_1$,$\beta_2$ and $\beta_3$ are small then regular fractions are preferred (see also Table~\ref{tab:23}). However, when at least one $|\beta_i|$ is large,  the regular fractions may not be optimal.
}\end{example}

\begin{figure}
\centering
\makebox{\includegraphics[scale=.43,angle=0]{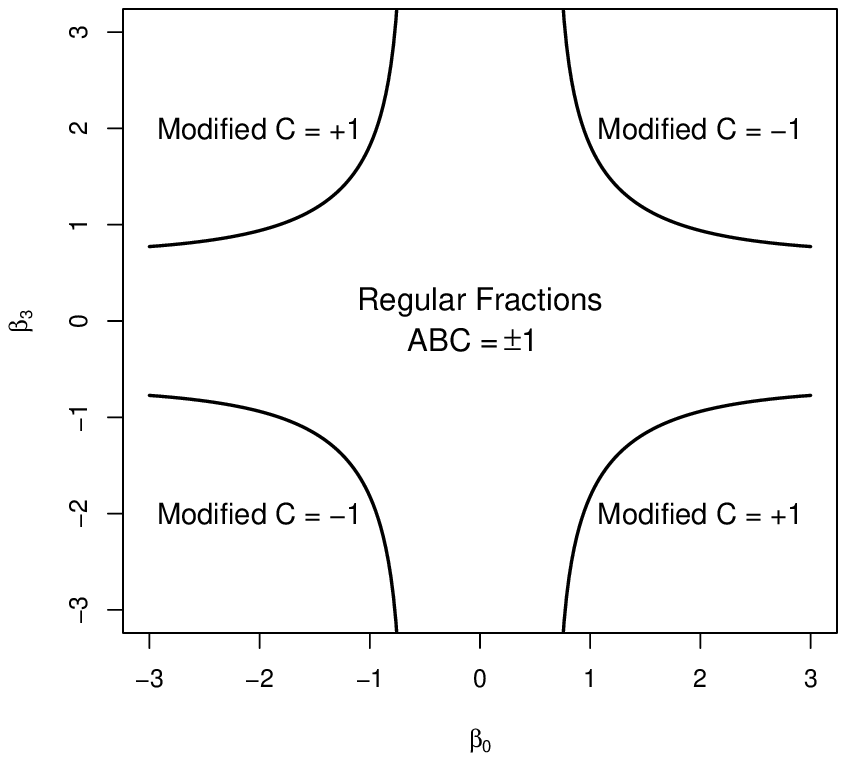}}\makebox{\includegraphics[scale=.43,angle=0]{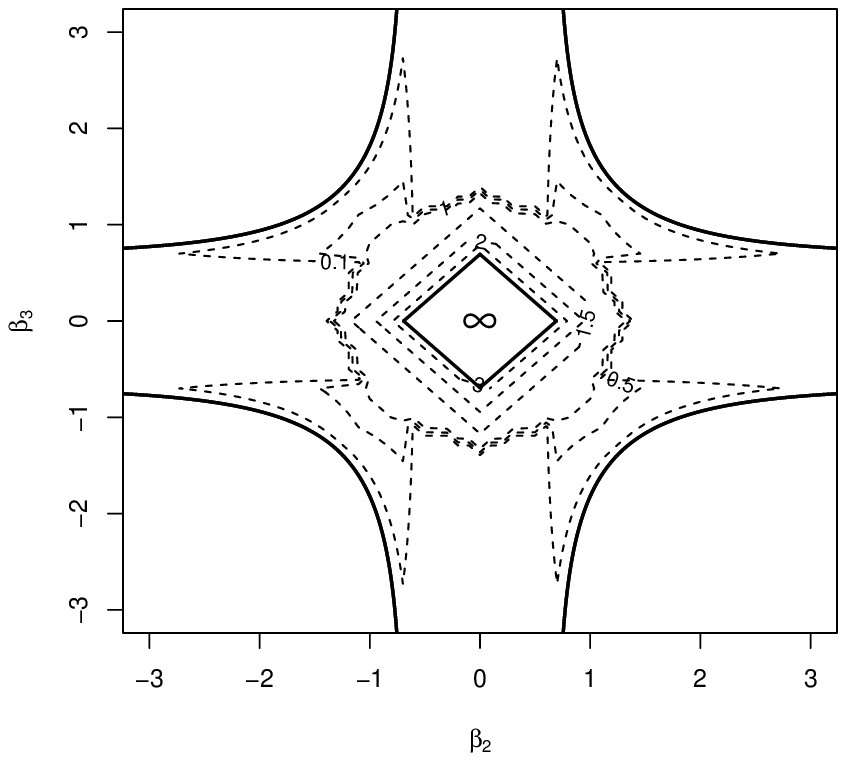}}
\caption{\label{fig0}Partitioning of the parameter space}
\end{figure}

In general, when all the $\beta_i$'s are nonzero, the regular fractions given by the rows $\{1,4,6,7\}$ or $\{2,3,5,8\}$ are not necessarily the optimal half-fractions. To explore this, we simulate the regression coefficients $\beta_0,$ $\beta_1,$ $\beta_2, \beta_3$ independently from different distributions and calculate the corresponding ${\mathbf w}$'s under logit, probit and complementary log-log links 10,000 times each. For each ${\mathbf w}$, we find the best (according to D-criterion) design supported on 4 distinct rows of the {model matrix}. By Lemma~\ref{xwx}, any such design has to be uniform. Table~\ref{tab:23} gives the percentages of times each of those designs turn out to be the optimal ones for the logit model (the results are somewhat similar for the other links).  It shows that the regular fractions are optimal when the $\beta_i$'s are close to zero.  In Table~\ref{tab:23}, we only report the non-regular fractions which turn out to be D-optimal for more than 15\% of the times. For the $2^4$ case, the results are similar, that is, when the $\beta_i$'s are nonzeros, the performance of the regular fractions given by $1=\pm ABCD$ are not very efficient in general.

We have done a simulation study to determine the efficiency of fractions, especially the regular ones. In order to describe a measure of efficiency, let us denote the D-criterion value as $\psi({\mathbf p}, {\mathbf w}) = |X'WX|$ for given ${\mathbf w} = (w_1, \ldots, w_{2^k})'$ and ${\mathbf p} = (p_1, \ldots, p_{2^k})'$. Suppose ${\mathbf p}_w$ is a D-optimal allocation with respect to ${\mathbf w}$. Then the loss of efficiency of ${\mathbf p}$ (with respect to a D-optimal allocation ${\mathbf p}_w$) given ${\mathbf w}$ can be defined as
\begin{eqnarray}\label{relloss}
R({\mathbf p}, {\mathbf w}) &=& 1 - \left(\frac{\psi({\mathbf p}, {\mathbf w})}{\psi({\mathbf p}_w, {\mathbf w})}\right)^{\frac{1}{d+1}}.
\end{eqnarray}

In Table~\ref{tab:23}, we provide within parentheses (the first number) the percentages of times that the regular fractions are at least 70\% efficient compared to the best half-fractions (it would correspond to the case where 42\% more runs are needed due to a poor choice of design). The second number within the parentheses is the median efficiency. It is clear that when the regular fractions are not D-optimal, they are usually not highly efficient either.

\begin{remark}\label{remarkpriorsimul0}
{\rm For each of the five situations described in Table~\ref{tab:23}, we also calculate the corresponding EW D-optimal half-fractions. For all five cases including the highly asymmetric fifth scenario, the regular fractions are EW D-optimal half-fractions.}\end{remark}

\begin{table}\caption{\label{tab:23}Distribution of D-optimal half-fractions  under $2^3$ main-effects model}
\centering
\fbox{%
{\footnotesize
\begin{tabular}{c|rcccc|c}
Rows &&\multicolumn{5}{c}{Percentages} \\
\hline
& $\beta_0\sim$ & \multicolumn{4}{c|}{$U(-10,10)$}       & \multicolumn{1}{c}{$N(0,5)$}      \\
Simulation & $\beta_1\sim$ & $U(-.3,.3)$  &   $U(-3,3)$  & $U(-3,0) $  &  $U(0,1)$ & $N(1,1)$ \\
Setup & $\beta_2\sim$      & $U(-.3,.3)$  &   $U(0,3) $  & $U(0,3)  $  &  $U(0,3)$ & $N(2,1)$  \\
& $\beta_3\sim$            & $U(-.3,.3)$  &   $U(1,5) $  & $U(-2,2) $  &  $U(0,5)$ & $N(3,1)$  \\
\hline
% 1467 & & 42.65  (99.5) &         0.07 (0.22)      &   0.86  (2.18)       &    0.95 (2.61) & 0.04 (2.89)  \\
% 2358 & & 42.02  (99.5) &         0.04 (0.23)      &   0.68  (2.12)       &    1.04 (2.56) & 0.08 (2.86)  \\
\multirow{2}{*}{1467}  & & 47.89          &         0.07            &   0.86        &    0.95       & 0.04  \\
                       & &(100,99.9)      &         (1.6,15.0)      &   (8.7,29.2)  &    (8.8,25.9) & (1.7,18.7)\\
\multirow{2}{*}{ 2358} & & 42.02          &         0.04            &   0.68        &    1.04       & 0.08  \\
                       & &(100,99.9)      &         (1.6,15.2)      &   (8.9,29.1)  &    (8.7,25.9) & (1.8,18.6)\\
\hline
 1235 & &  &      16.78  &          &    35.62              & 21.50     \\
 1347 & &  &             &   19.98  &                       &           \\
 1567 & &  &      17.45  &   19.21  &                       &           \\
 2348 & &  &      17.54  &   19.11  &                       &           \\
 2568 & &  &             &   20.01  &                       &           \\
 4678 & &  &      16.12  &          &    35.41              & 21.65     \\
\end{tabular}
}}
\end{table}

\begin{remark}\label{remarkpriorsimul1}
{\rm In Table~\ref{tab:cpu} and later (Table~\ref{tab:23} and Table~\ref{tab:odor2}) we have used distributions for $\boldsymbol{\beta}$ in two ways. For locally D-optimal designs these distributions are used to simulate the assumed values in order to study the properties of the designs, especially robustness. For EW D-optimal designs these distributions are used as priors.}\end{remark}

\begin{remark}\label{remarkpriorsimul2}
{\rm The priors for $\boldsymbol{\beta}$ should be chosen carefully for real applications. For example, a uniform prior on  $\beta_i\sim [-a, a]$ would indicate that the experimenter does not know much about the corresponding factor. If $\beta_i\sim [0, b]$ then the experimenter knows the direction of the corresponding factor effect.  In our odor study example, factor $A$ (algae) has two levels: raffinated or solvent extracted algae ($-1$) and catfish pond algae ($+1$). The scientists initially assessed that raffinated algae has residual lipid which should prevent absorber to interact with volatiles, causing odor to release. Hence it is expected that $\beta_i$ for this factor should be nonnegative. In this case, one may take the prior on $[0,b]$. On the other hand, for factor $B$ (Scavenger), it is not known before conducting the experiment whether Activated Carbon ($-1$) is better or worse than Zeolite ($+1$). In this case, a symmetric prior on $[-a, a]$ would be more appropriate.
}\end{remark}

\begin{remark}\label{remarkfraction}
{\rm Consider the problem of obtaining the locally D-optimal fractional factorial designs when the number of experimental settings ($m$, say) is fixed. If the total number of factors under consideration is not too large, one can always calculate the D-efficiencies of all fractions and choose the best one. However, this is a computationally expensive strategy for large $k$'s so we need an alternative. One such strategy would be to choose the $m$ largest $w_i$'s and the corresponding rows, since {those $w_i$ represent the information at the corresponding design points}. Another one would be to use our algorithms discussed in Section~3.3 to find an optimal allocation for the full factorial designs first, then to choose the $m$ largest $p_i$'s and scale them appropriately. One has to be careful, however, in order to avoid designs which would not allow the estimation of the model parameters. In this case, the exchange algorithm described in Section~\ref{algo2} may be used to choose the fraction with given $m$ experimental units. Our simulations (not presented here) show that both of these methods perform satisfactorily with the second method giving designs which are generally more than $95\%$ efficient for four factors with the main-effects model. This method will be used for computations in the next section.
}\end{remark}

\fontsize{10.95}{14pt plus.8pt minus .6pt}\selectfont
\setcounter{chapter}{5}
\setcounter{equation}{0} %-1
\setcounter{section}{0} %-1

{\bf 5. Robustness}

In this section, we will study the robustness of locally D-optimal designs over the assumed parameter values.

\setcounter{section}{1}

{\bf 5.1 Most robust minimally supported designs}\label{sec:robustsaturated}

Minimally supported designs have been studied extensively. For continuous or quantitative factors, these designs can be D-optimal for many linear and non-linear models. In our setup of qualitative factors, these designs are attractive since they use the minimal number, $d+1$, of experimental conditions. In many applications, fewer experimental conditions are desirable. In this section, we will examine the robustness of minimally supported designs. Our next result gives necessary and sufficient conditions for a fraction to be a D-optimal minimally supported design. Note that Theorem~\ref{mfractioncorollary} is an immediate consequence of Lemma~\ref{xwx}.

\begin{theorem}\label{mfractioncorollary}
Let ${\mathbf I}=\{i_1,\ldots,i_{d+1}\}\subset \{1,\ldots,2^k\}$ be an index set. A design ${\mathbf p}_I=(p_1,\ldots,$ $p_{2^k})^\prime$ satisfying $p_i=0, \forall i\notin I$ is $D$-optimal among minimally supported designs if and only if $$ p_{i_1}=\cdots=p_{i_{d+1}}=\frac{1}{d+1}\mbox{ and }{\mathbf I}\mbox{ maximizes } |X[i_1, \ldots, i_{d+1}]|^2w_{i_1}\cdots w_{i_{d+1}}~.$$
\end{theorem}

Recall that we denoted the loss of efficiency of ${\mathbf p}$ in (\ref{relloss}) by $R({\mathbf p}, {\mathbf w})$. For investigating the robustness of a design, let us define the maximum loss of efficiency of a given design ${\mathbf p}$ with respect to a specified region ${\cal W}$ of ${\mathbf w}$ by
\begin{equation}\label{rmaxdef}
R_{\max}({\mathbf p}) = \max_{{\mathbf w} \in {\cal W}} R({\mathbf p}, {\mathbf w}).
\end{equation}

It can be shown that the region ${\cal W}$ takes the form of $[a,\ b]^{2^k}$ for $2^k$ main-effects model if the range of each of the regression coefficients is an interval symmetric about 0. For example, for a $2^4$ main-effects model, if all the regression coefficients range between $[-3,3]$, then ${\cal W} =[3.06\times 10^{-7},\ 0.25]^{16}$ for logit link, and $[8.33\times 10^{-49}, 0.637]^{16}$ for probit link. This is the rationale for the choice of the range of $w_i$'s in Theorem~\ref{drobusttheorem} below. A design which minimizes the maximum loss of efficiency will be called {\it most robust}. {Note that this criterion is also known as ``maximin efficiency'' in the literature (see, for example, Dette (1997)). For unbounded $\beta_i$'s with a prior distribution, one may use $.99$ or $.95$ quantile instead of the maximum loss to measure the robustness.}

\begin{theorem}\label{drobusttheorem}
Suppose $k\geq 3$ and $\ d(d+1)\leq 2^{k+1}-4$. Suppose $w_i \in [a,\ b]$, $i=1, \ldots, 2^k$, $0 < a < b$. Let ${\mathbf I}=\{i_1,\ldots,i_{d+1}\}$ be an index set which maximizes $|X[i_1, i_2, \ldots,i_{d+1}]|^2$. Then the design ${\mathbf p}_I=(p_1,\ldots,p_{2^k})'$ satisfying $p_{i_1}=\cdots=p_{i_{d+1}}=\frac{1}{d+1}$ is a most robust minimally supported design with
maximum loss $1-\frac{a}{b}$ in efficiency compared to other minimally supported designs.
\end{theorem}

Based on Theorem~\ref{drobusttheorem}, the maximum loss of efficiency depends on the range of $w_i$'s. The result is meaningful only if the interval $[a,b]$ is bounded away from 0. Figure~\ref{fig1} provides some idea about the possible bounds of $w_i$'s for commonly used link functions.  For example, for $2^3$ designs with main-effects model, if $0.105 \le w_i \le 0.25$ under logit link (see Remark 4.1.1 of Yang et al.~(2012)), then the maximum loss of efficiency of the regular half-fractional design satisfying $p_1 = p_4 = p_6 = p_7 = 1/4$ is $1 - 0.105/0.25 = 58\%$. The more certain we are about the range of $w_i$'s, the more useful the result will be.

Note that for $k=2$ all $4$ minimally supported designs perform equally well (or equally badly). So they are all most robust under the above definition. For main-effects models, the condition $d(d+1)\leq 2^{k+1}-4$ in Theorem~\ref{drobusttheorem} is guaranteed whenever $k\geq 3$. A most robust minimally supported design can be obtained by searching for an index set $\{i_1,\ldots,i_{d+1}\}$ which maximizes $|X[i_1, i_2, \ldots,i_{d+1}]|^2$. Note that such an index set is usually not unique. Based on Lemma~\ref{indexsetlemma}, if the index set $\{i_1, \ldots, i_{d+1}\}$ maximizes $|X[i_1, \ldots, i_{d+1}]|^2$, then there always exists another index set $\{i_1', \ldots, i_{d+1}'\}$ such that $|X[i_1, \ldots, i_{d+1}]|^2 = |X[i_1', \ldots, i_{d+1}']|^2$. It should also be noted that a most robust minimally supported design may involve a set of experimental conditions $\{i_1,\ldots,i_{d+1}\}$ which does not maximize $|X[i_1,\ldots,i_{d+1}]|^2$. For example, consider a $2^{3-1}$ design with main-effects model. Suppose $w_i\in [a,b]$, $i=1,\ldots,8$. If $4a>b$, then the most robust minimally supported designs are the $2^{3-1}$ regular fractions. Otherwise,  if $4a\leq b$, then any uniform design restricted to $\{i_1,i_2,i_3,i_4\}$ satisfying $|X[i_1,i_2,i_3,i_4]|$ $\neq 0$ is a most robust minimally supported design.

%\clearpage
\setcounter{section}{2}

  {\bf 5.2 Robustness of uniform designs}\label{sec:robustuniform}

As mentioned before, for examples in Section~3.2, a design is called ``uniform'' if the allocation of experimental units is the same for all points in the support of the design. Yang, Mandal and Majumdar (2012) showed that for a $2^2$ main-effects model, the uniform design is the most robust design in terms of maximum loss of efficiency. In this section, we use simulation studies to examine the robustness of uniform designs and EW D-optimal designs for higher order cases.

For illustration, we use a $2^4$ main-effects model. We simulate $\beta_0, \ldots, \beta_4$ from different distributions 1000 times each and calculate the corresponding ${\mathbf w}$'s, denoted by vectors ${\mathbf w}_1$, $\ldots$, ${\mathbf w}_{1000}$~. For each ${\mathbf w}_s$, we use the algorithm described in Section~\ref{sec:robustsaturated} to obtain a D-optimal allocation ${\mathbf p}_s$~. For any allocation ${\mathbf p}$, let  $R_{100\alpha}({\mathbf p})$ denote the $\alpha$th quantile of the set of loss of efficiencies $\{R({\mathbf p}, {\mathbf w}_s), \ s=1,\ldots, 1000\}$. Thus $R_{100}({\mathbf p})=R_{\max}({\mathbf p})$ which is the $R_{\max}$ defined in (\ref{rmaxdef}) with ${\cal W} = \{{\mathbf w}_1, \ldots, {\mathbf w}_{1000}\}$. The quantities $R_{99}({\mathbf p})$ and $R_{95}({\mathbf p})$ are more reliable in measuring the robustness of ${\mathbf p}$.

\begin{table}\caption{\label{tab:robust}Loss of efficiency of $2^4$ uniform design}
\centering
\fbox{%
{\scriptsize
\begin{tabular}{c|rccccccccccc}
 & \multicolumn{12}{c}{Percentages} \\
           & \multicolumn{3}{c}{$\beta_0\sim U(-3,3)$}  & \multicolumn{3}{c}{$U(-1,1)$} & \multicolumn{3}{c}{$U(-3,0)$}      & \multicolumn{3}{c}{$N(0,5)$}  \\
           & \multicolumn{3}{c}{$\beta_1\sim U(-1,1)$}  & \multicolumn{3}{c}{$U(0,1)$}  & \multicolumn{3}{c}{$U(1,3)$}      & \multicolumn{3}{c}{$N(0,1)$}  \\
Simulation & \multicolumn{3}{c}{$\beta_2\sim U(-1,1)$}  & \multicolumn{3}{c}{$U(0,1)$}  & \multicolumn{3}{c}{$U(1,3)$}      & \multicolumn{3}{c}{$N(2,1)$}  \\
Setup      & \multicolumn{3}{c}{$\beta_3\sim U(-1,1)$}  & \multicolumn{3}{c}{$U(0,1)$}  & \multicolumn{3}{c}{$U(-3,-1)$}      & \multicolumn{3}{c}{$N(-.5,2)$}\\
           & \multicolumn{3}{c}{$\beta_4\sim U(-1,1)$}  & \multicolumn{3}{c}{$U(0,1)$}  & \multicolumn{3}{c}{$U(-3,-1)$}      & \multicolumn{3}{c}{$N(-.5,2)$}\\
\hline
Quantiles & (I) & (II) & (III)                     & (I) & (II) & (III)     & (I) & (II) & (III)      & (I) & (II) & (III) \\
 $R_{99} $    &  .348 & .353 &  .348  & .146 & .111 & .112  & .503 &  .273 & .299  &  .650  & .864 & .726 \\
 $R_{95} $    &  .299 & .304 &  .299  & .128 & .094 & .093  & .495 &  .251 & .256  &  .617  & .788 & .670 \\
 $R_{90} $    &  .271 & .274 &  .271  & .117 & .084 & .085  & .488 &  .239 & .233  &  .589  & .739 & .629 \\
\hline
 \multicolumn{13}{c} {Note: (I) = $R_{100\alpha}({\mathbf p}_u)$, (II) = $\displaystyle{\min_{1\leq s\leq 1000}} R_{100\alpha}({\mathbf p}_s)$, (III) = $R_{100\alpha}({\mathbf p}_e)$.}\\
 \multicolumn{13}{c} { {${\mathbf p}_u$ is the uniform design, ${\mathbf p}_s$ is the locally D-optimal design and ${\mathbf p}_e$ is the EW D-optimal design.}}
\end{tabular}
}}\hspace{11in}
\end{table}

Table~\ref{tab:robust} compares the $R_{100\alpha}$ of the uniform design ${\mathbf p}_u = (1/16,\ldots,1/16)^\prime$ with the minimum of $R_{100\alpha}({\mathbf p_s})$ for the optimal allocations ${\mathbf p}_s$, $s=1, \ldots, 1000$, as well as the $R_{100\alpha}$ of the EW design ${\mathbf p}_e$~. In this table, if the values of column (I) is smaller than those of column (II), then we can conclude that the uniform design is better than all the D-optimal designs in terms of the quantiles of loss of efficiency. This happens in many situations. Table~\ref{tab:robust} provides strong evidence for fact that the uniform design ${\mathbf p}_u$ is one of the most robust ones if the $\beta_i$'s are expected to come from an interval that is symmetric around zero. This is consistent with the conclusion of Cox (1988).

However, there are situations where the uniform design does not perform well, as illustrated by the two middle blocks of Table~\ref{tab:robust}. If the signs of the regression coefficients are known, it is advisable not to use the uniform design. For many practical applications, the experimenter will have some idea of the direction of effects of factors, which in statistical terms determines the signs of the regression coefficients. For these situations, it turns out that the performance of the EW D-optimal designs is comparable to that of the most robust designs, even when the uniform design does not perform well (see columns (III) in Table~\ref{tab:robust}, where ${\mathbf p}_e$ is the EW design). Hence we recommend the use of EW D-optimal designs when the experimenter has some idea about the signs of $\beta_i$'s. Uniform designs are recommended in the absence of prior knowledge of the sign of the regression parameters.

Now consider the uniform designs restricted to regular fractions. Again we use $2^4$ main-effects model as illustration and consider the uniform designs restricted to the regular half-fractions identified by $1=\pm ABCD$. We performed simulations as above and our conclusions are similar, that is, uniform designs on regular fractions are among the most robust ones if the signs of the regression parameters are unknown but they may not perform well if the signs of $\beta_i$'s are known.

%\clearpage

\fontsize{10.95}{14pt plus.8pt minus .6pt}\selectfont
\setcounter{chapter}{6}
\setcounter{equation}{0} %-1
\setcounter{example}{0} %-1
\setcounter{section}{1}

  {\bf 6. Examples}

In this section we discuss two examples.

\begin{example}{\rm First we revisit the odor examples discussed in the introduction. The $2^{4-1}_{IV}$ design given by $D=-ABC$ was used with 5 replications per experimental setup. For factor $C$, the polypropylene used in this experiment is in tiny crystal form as opposed to fine powder which leads the scientist to speculate that $\beta_3$ should be positive. Moreover one expects that the presence of compatabilizers should reduce the odor and hence $\beta_4$ is expected to be positive. Initial results from the experiment indicate that the number of successes is increasing in the level of $A$ (from $-1$ to $+1$). Let us examine the efficiency of the design used in this experiment in view of these facts and consider an EW D-optimal design with the following ranges, ($-3,3$) for $\beta_0,\beta_2$ and (0,3) for $\beta_1, \beta_3, \beta_4$. Note that these priors are reasonably uninformative except for the directions of effects of the factors (signs of the parameters). Furthermore, if the design points are not restricted to the original half-fraction, the best EW D-optimal design with 40 experimental units, given by ${\mathbf n}_{EW}$, is supported on 13 points.

\begin{table}[ht]
\caption{\label{tab:odor}Optimal design for the Odor Study}
\begin{center}
\begin{tabular}{cccccccc}
    A & B & C & D & E($w_i$) &  ${\mathbf n}_{odor}$ & ${\mathbf n}_{EW}$ & ${\mathbf n}_{EW{\frac{1}{2}}}$ \\
\hline
  $+1$& $+1$&   $+1$&  $+1$&  0.050  &         &        &         \\
  $+1$& $+1$&   $+1$&  $-1$&  0.105  &     5   &   3    &     7   \\
  $+1$& $+1$&   $-1$&  $+1$&  0.105  &     5   &   4    &     3   \\
  $+1$& $+1$&   $-1$&  $-1$&  0.105  &         &   3    &         \\
  $+1$& $-1$&   $+1$&  $+1$&  0.050  &     5   &        &         \\
  $+1$& $-1$&   $+1$&  $-1$&  0.105  &         &   4    &         \\
  $+1$& $-1$&   $-1$&  $+1$&  0.105  &         &   3    &     4   \\
  $+1$& $-1$&   $-1$&  $-1$&  0.105  &     5   &   3    &     6   \\
  $-1$& $+1$&   $+1$&  $+1$&  0.105  &     5   &   4    &         \\
  $-1$& $+1$&   $+1$&  $-1$&  0.105  &         &   3    &     3   \\
  $-1$& $+1$&   $-1$&  $+1$&  0.105  &         &   2    &     7   \\
  $-1$& $+1$&   $-1$&  $-1$&  0.050  &     5   &   1    &         \\
  $-1$& $-1$&   $+1$&  $+1$&  0.105  &         &   3    &     6   \\
  $-1$& $-1$&   $+1$&  $-1$&  0.105  &     5   &   3    &     4   \\
  $-1$& $-1$&   $-1$&  $+1$&  0.105  &     5   &   4    &         \\
  $-1$& $-1$&   $-1$&  $-1$&  0.050  &         &        &         \\
\hline
\end{tabular}
\end{center}
\end{table}

In order to compare the performance of the three designs given in Table~\ref{tab:odor}, we draw 1000 random samples of the $\beta_i$'s from the setup discussed above and for each of them calculate the locally D-optimal design with 40 runs. Then we calculate the loss of efficiencies of the EW D-optimal design (${\mathbf n}_{EW}$) and EW D-optimal half-fraction (${\mathbf n}_{EW{\frac{1}{2}}}$) as well as that of the original design used (${\mathbf n}_{odor}$), with respect to the locally D-optimal design. The mean, standard deviation and some quantiles of the loss of efficiencies are given in Table~\ref{tab:odor2}. These numbers indicate that the EW D-optimal design is around 20\% more efficient than the original one, while the EW half-fraction design is about 10\% more efficient than the original one.

\begin{table}[ht]
\caption{\label{tab:odor2}Odor Study: Loss of efficiencies of different designs}
\begin{center}
\begin{tabular}{lccccr}
             Design                               &  $R_{99}$   &  $R_{95}$     &  $R_{90}$   &  Mean & SD \\
             \hline
EW design (${\mathbf n}_{EW}$)                    & 51.4 & 46.6 & 44.7 & 33.0  &  9.5   \\
EW half-fraction (${\mathbf n}_{EW{\frac{1}{2}}}$)& 77.2 & 69.5 & 63.2 & 41.9  & 15.7   \\
Original design (${\mathbf n}_{odor}$)            & 84.8 & 76.8 & 70.1 & 51.8  & 15.1   \\
\hline
\end{tabular}
\end{center}
\end{table}

}\end{example}

\begin{example}{\rm
Hamada and Nelder (1997) discussed a $2^{4-1}$ fractional factorial experiment performed at IIT Thompson laboratory that was originally reported by Martin, Parker and Zenick (1987). This was a windshield molding slugging experiment where the outcome was whether the molding was good or not. There were four factors each at two levels: ($A$) poly-film thickness (0.0025, 0.00175), ($B$) oil mixture ratio (1:20, 1:10), ($C$) material of gloves (cotton, nylon), and ($D$) the condition of metal blanks (dry underside, oily underside). By analyzing the data presented in Hamada and Nelder (1997), we get an estimate of the unknown parameter as  $\hat{\boldsymbol\beta}=(1.77,-1.57,0.13,-0.80,-0.14)^\prime$ under logit link. If one wants to conduct a follow-up experiment on half-fractions, then it is sensible to use the knowledge obtained by analyzing the data. With the knowledge of  $\hat{\boldsymbol\beta}$, let us take the assumed value of ${\boldsymbol\beta}$ as $(2,-1.5,0.1,-1,-0.1)^\prime$. The locally D-optimal design ${\mathbf p}_a$ is given in  Table~\ref{tab:hn}. Another option is to consider a range for the possible values of the regression parameters, namely, $(1,3)$ for $\beta_0$, $(-3,-1)$ for $\beta_1$, $(-0.5,0.5)$ for $\beta_2, \beta_4$, and $(-1,0)$ for $\beta_3$. For this choice of range for the parameter values with independence and uniform distributions, the EW D-optimal half-fractional design ${\mathbf p}_e$ is also given in Table~\ref{tab:hn}. We have calculated the linear predictor $\eta$ and success probability $\pi$ for all possible experimental settings. It seems that a good fraction would not favor high success probabilities very much. This is one of the main differences between the design reported by  Hamada and Nelder (denoted by ${\mathbf p}_{HN}$) and our designs (denoted by ${\mathbf p}_a$ and ${\mathbf p}_e$). Note that these two designs have six rows in common. The last two columns of Table~\ref{tab:hn} give the Baysian D-optimal and EW D-optimal designs, respectively. It can be seen that the optimal allocation for these two designs are quite similar, and both of them are supported on the same rows.

\begin{table}
\caption{\label{tab:hn}Optimal half-fraction design for Windshield Molding Experiment}
\centering
\fbox{%
\begin{tabular}{cccccrcccccc}
Row & A & B & C & D & $\eta$ & $\pi$ & ${\mathbf p}_{HN}$ & ${\mathbf p}_a$ & ${\mathbf p}_e$ & ${\mathbf p}_B$ & ${\mathbf p}_{e_f}$ \\
\hline
5  & $+1$  & $-1$  & $+1$  & $+1$  & -0.87  & 0.295  &         &   0.044  & 0.184 &  0.073   &     0.092  \\
1  & $+1$  & $+1$  & $+1$  & $+1$  & -0.61  & 0.352  &  0.125  &   0.178  & 0.011 &  0.117   &     0.103  \\
6  & $+1$  & $-1$  & $+1$  & $-1$  & -0.59  & 0.357  &  0.125  &   0.178  & 0.011 &  0.118   &     0.103  \\
2  & $+1$  & $+1$  & $+1$  & $-1$  & -0.33  & 0.418  &         &   0.059  & 0.184 &  0.078   &     0.092  \\
7  & $+1$  & $-1$  & $-1$  & $+1$  &  0.73  & 0.675  &  0.125  &   0.163  &       &  0.125   &     0.103  \\
3  & $+1$  & $+1$  & $-1$  & $+1$  &  0.99  & 0.729  &         &          & 0.195 &  0.079   &     0.091  \\
8  & $+1$  & $-1$  & $-1$  & $-1$  &  1.01  & 0.733  &         &          & 0.195 &  0.078   &     0.091  \\
4  & $+1$  & $+1$  & $-1$  & $-1$  &  1.27  & 0.781  &  0.125  &   0.147  &       &  0.115   &     0.103  \\
13 & $-1$  & $-1$  & $+1$  & $+1$  &  2.27  & 0.906  &  0.125  &   0.158  & 0.111 &  0.061   &     0.054  \\
9  & $-1$  & $+1$  & $+1$  & $+1$  &  2.53  & 0.926  &         &          &       &  0.053   &     0.057  \\
14 & $-1$  & $-1$  & $+1$  & $-1$  &  2.55  & 0.928  &         &          &       &  0.043   &     0.057  \\
10 & $-1$  & $+1$  & $+1$  & $-1$  &  2.81  & 0.943  &  0.125  &   0.074  & 0.110 &  0.061   &     0.053  \\
15 & $-1$  & $-1$  & $-1$  & $+1$  &  3.87  & 0.980  &         &          &       &          &            \\
11 & $-1$  & $+1$  & $-1$  & $+1$  &  4.13  & 0.984  &  0.125  &          &       &          &            \\
16 & $-1$  & $-1$  & $-1$  & $-1$  &  4.15  & 0.984  &  0.125  &          &       &          &            \\
12 & $-1$  & $+1$  & $-1$  & $-1$  &  4.41  & 0.988  &         &          &       &          &            \\
\end{tabular}
}

\vspace{.1in}

 {Notation:  ${\mathbf p}_{HN}$: Design reported by Hamada and Nelder, ${\mathbf p}_a$: Locally D-optimal design, ${\mathbf p}_e$: EW D-optimal half-fraction, ${\mathbf p}_B$: Bayesian D-optimal design, ${\mathbf p}_{e_f}$: EW D-optimal design}
\end{table}
}\end{example}

\fontsize{10.95}{14pt plus.8pt minus .6pt}\selectfont
\setcounter{chapter}{7}
\setcounter{equation}{0} %-1

  {\bf 7. Discussion and Future Research}

For binary response, the logit link is the most commonly used link in practice. The situation under this link function is close to that in the linear model case because typically $w_i$'s are not too close to $0$ and do not vary much. Similar to the cases of linear models, uniform designs perform well under logit link, more than other popular link functions. In general, the performance of the logit and probit links are similar, while that of the complementary log-log link is somewhat different from others. For example, if we consider a $2^2$ experiment with a main-effects model, the efficiency of the uniform design with respect to the Bayes D-optimal design is 99.99\% under logit link, but is only 89.6\% under complementary log-log link. Figure~\ref{fig1} provides a graphical display of the weight function $(w)$ for commonly used link functions. As seen from the figure, complementary log-log link function is not symmetric about 0. This partly explains the poor performance of the uniform design under this link. Nevertheless, the EW D-optimal designs are still highly efficient across different link functions. For the same setup, the efficiencies of EW designs with respect to the corresponding Bayesian D-optimal designs are $99.99\%$ (logit link), $99.94\%$ (probit link), $99.77\%$ (log-log link), and $100.00\%$ (complementary log-log link), respectively. From all of our simulations it appears that EW D-optimal designs are excellent surrogates of Bayes D-optimal designs. A more extensive investigation is planned for the future.

\begin{figure}[h]
\centering
\makebox{\includegraphics[scale=.4,angle=0]{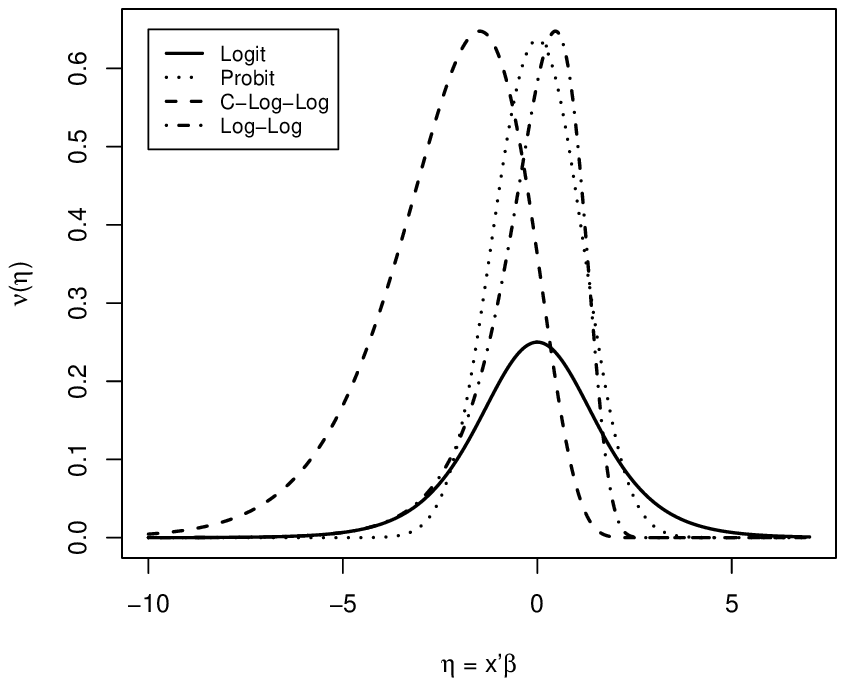}}
\caption{\label{fig1}$w_i = \nu(\eta_i) = \nu({\mathbf x}_i'\boldsymbol\beta)$ for commonly used link functions}
\end{figure}

 {It should also be noted that the efficiencies depend on the priors used for the parameters, and hence the prior on the $\boldsymbol\beta$s should  be different for different link functions in order to maintain roughly consistent prior beliefs about the success probabilities under different experimental setups.}

Our recommendation is to use EW D-optimal designs unless the experimenter has absolutely no prior knowledge of the parameters, in which case it is recommended to use the uniform design. In EW optimality, we replace the $w_i$'s by their expectations. It may be noted, however, that taking the average of $w_i$'s is not same as taking the average of $\beta_i$'s. Let us illustrate this with a $2^4$ design with main-effects model. Table~\ref{tab5} below uses the notations from Table~\ref{tab:robust}. Suppose $\beta_0 \sim U(-3,0)$, $\beta_1, \beta_2 \sim U(1,3)$, $\beta_3, \beta_4 \sim U(-3,-1)$, and the $\beta_i$'s are independent. It is clear that the uniform design performs much worse compared to the most robust design, while the performance of the EW D-optimal design is comparable with the best design. The last column corresponds to the locally D-optimal design where the assumed value of the parameter is taken to be the midpoints of the ranges of $\beta_i$'s mentioned above. Clearly this is worse than the EW D-optimal design.

\begin{table}[h]\caption{\label{tab5}Loss of efficiencies of different designs for $2^4$ main-effects model}
\centering
\fbox{%
{\begin{tabular}{ccccc}
          & Uniform & Most robust & EW D-opt & $E(\beta)$ D-opt\\
\hline
%$R_{100}$ & 0.507 & 0.287 & 0.316 &   0.350  \\
 $R_{99}$  & 0.503 & 0.273 & 0.299 &   0.331  \\
 $R_{95}$  & 0.495 & 0.251 & 0.256 &   0.284  \\
 $R_{90}$  & 0.488 & 0.239 & 0.233 &   0.251  \\
%$R_{75}$  & 0.472 & 0.210 & 0.204 &   0.210  \\
%$R_{50}$  & 0.447 & 0.171 & 0.169 &   0.167  \\
\end{tabular}
}}\end{table}

In the linear model setup, as the potential columns in the model matrix are orthogonal, analysis of experimental data based on regular fractions is not unduly biased by the omission of non-negligible model terms. Under a GLM setup, the regular fractions may give larger than necessary variance for some models. In this paper, we did not consider the performance of different designs under model robustness. Moreover, because of the bias-variance trade-off, regular fractions (or other designs) may not be model-robust. Extending optimal designs based on GLMs to topics such as confounding, aberration, and trade-off between variance and bias represents an important topic for future research.

%\clearpage

 \medskip {\bf Acknowledgment:}

We thank Dr. Suraj Sharma for providing the details of the Odor Study, and Dr. Yaming Yu for sharing codes for the Multiplicative and Cocktail algorithms. We also thank the reviewers for comments and suggestions that substantially improved the quality of the manuscript. This research is in part supported by NSF Grant DMS-09-05731 and NSA Grant H98230-13-1-0251. \par

%\clearpage

%\renewcommand{\theequation}{A.\arabic{equation}}
%\appendix

%\begin{center}\noindent {\large\bf Appendix}\end{center}

%\noindent {\textbf{\emph{1. Proof of Theorem~\ref{thm:1}}}}

\singlespacing

\section*{References}
\begin{description}

\item Agresti, A. (2002). {\em Categorical Data Analysis}, Second Edition. John Wiley \& Sons, New York.

\item Atkinson, A. C., Donev, A. N. and Tobias, R. D. (2007). {\em Optimum Experimental Designs, with SAS}, Oxford University Press.

\item Chaloner, K. and Verdinelli, I. (1995). ``Bayesian experimental design: a review", {\em Statistical Science} {\bf 10}, 273$-$304.

\item Chernoff, H. (1953). ``Locally optimal designs for estimating parameters", {\em Annals of Mathematical Statistics} {\bf 24}, 586$-$602.

\item Cox, D. R. (1988). ``A Note on Design when Response has an Exponential Family Distribution", {\em Biometrika}, {\bf 75}, 161$-$164.

\item Dette, H. (1997). ``Designing experiments with respect to `standardized' optimality criteria", {\em Journal of the Royal Statistical Society, Series B}, {\bf 59}, 97$-$110.

\item Dobson, A. J. and Barnett, A. (2008). {\em An Introduction to Generalized Linear Models}, Third Edition. Chapman and Hall/CRC, London.

\item Dorta-Guerra, R., Gonz\'{a}lez-D\'{a}vila, E. and Ginebra, J. (2008). ``Two-level experiments for binary response data'', {\em Computational Statistics and Data Analysis}, {\bf 53}, 196$-$208.

\item Dror, H. A. and Steinberg, D. M. (2006). ``Robust Experimental Design for Multivariate Generalized Linear Models'', {\em Technometrics}, {\bf 48}, 520$-$529.

\item Dror, H. A. and Steinberg, D. M. (2008).  ``Sequential Experimental Designs for Generalized Linear Models'', \emph{Journal of the American Statistical Association}, {\bf 103}, 288$-$298.

\item Fedorov, V.V. (1972). {\em Theory of Optimal Experiments}, Academic Press, New~York.

\item Fedorov, V.V. and Hackl, P. (1997). \emph{Model-Oriented Design of Experiments}, Springer, New~York.

\item Fedorov, V.V. and Leonov, S.L. (2014). {\em Optimal Design for Nonlinear Response Models}, Chapman \& Hall/CRC.

\item Ford, I., Titterington, D. M. and Kitsos, C. P. (1989). Recent advances in nonlinear experimental design, {\it Technometrics}, {\bf 31}, 49$-$60.

\item Gonz\'{a}lez-D\'{a}vila, E., Dorta-Guerra, R. and Ginebra, J. (2007), ``On the information in two-level experiments", {\em Model Assisted Statistics and Applications}, {\bf 2}, 173$-$187.

\item Gotwalt. C. M., Jones, B. A. and Steinberg, D. M. (2009), ``Fast Computation of Designs Robust to Parameter Uncertainty for Nonlinear Settings'',   {\em Technometrics},    {\bf 51}, 88$-$95.

\item Gra{\ss}hoff, U. and Schwabe, R. (2008). ``Optimal design for the Bradley-Terry paired comparison model", \emph{Statistical Methods and Applications},  {\bf 17}, 275$-$289.

%\item Grimshaw, S. D.,  Collings, B. J., Larsen, W. A. and Hurt C. R. (2001). ``Eliciting Factor Importance in a Designed Experiment",  \emph{Technometrics}, {\bf 43}, 133$-$146.

\item Hamada, M. and Nelder, J. A. (1997). ``Generalized linear models for quality-improvement experiments'', \emph{Journal of Quality Technology}, {\bf 29}, 292$-$304.

\item Imhof, L. A. (2001). ``Maximin designs for exponential growth models and heteroscedastic polynomial models'', {\em Annals of Statistics}, {\bf 29}, 561$-$576.

\item Kiefer, J. (1974). ``General equivalence theory for optimum designs (approximate theory)", {\em Annals of Statistics}, {\bf 2}, 849$-$879.

\item Khuri, A. I., Mukherjee, B., Sinha, B. K. and Ghosh, M. (2006). ``Design issues for generalized linear models: A Review", {\em Statistical Science}, {\bf 21}, 376$-$399.

\item Li, G. and Majumdar, D. (2008). ``D-optimal designs for logistic models with three and four parameters'', \emph{Journal of Statistical Planning and Inference}, {\bf 138}, 1950$-$1959.

\item Li, G. and Majumdar, D. (2009). ``Some results on D-optimal designs for nonlinear models with applications'', \emph{Biometrika}, {\bf 96}, 487$-$493.

\item Lindsey, J. (1997). {\em Applying Generalized Linear Models}. Springer, New~York.

\item Mandal, A., Wong, W. K. and Yu, Y. (2014). ``Algorithmic Searches for Optimal Designs'', {\em Handbooks on Modern Statistical Methods}, Chapman and Hall/CRC.

\item Martin, B., Parker, D. and Zenick, L. (1987). ``Minimize slugging by optimizing controllable factors on topaz windshield molding'', In: \emph{Fifth Symposium on Taguchi Methods, American Supplier Institute, Inc., Dearborn, MI}, 519$-$526.

\item McCullagh, P. and Nelder, J. (1989). {\em Generalized Linear Models},  Second Edition. Chapman and Hall/CRC, Boca Raton.

\item McCulloch, C., and Searle, S. (2001). {\em Generalized, linear and mixed models}. Wiley, New~York.

%\item Miller, A. and Sitter, R. R. (2001). ``Using the Folded-Over 12-Run Plackett-Burman Design to Consider Interactions", \emph{Technometrics}, {\bf 43}, 44$-$55.

\item M\"uller, W. G. (2007). {\em Collecting Spatial Data: Optimum Design of Experiments for Random Fields}, 3rd Edition, Springer.

\item Myers, R. M., Montgomery, D. C., and Vining, G. G. (2002). {\em Generalized Linear Models with Applications in Engineering and Statistics}. John Wiley, New~York.

\item Nair, V., Strecher, V., Fagerlin, A., Ubel, P., Resnicow, K., Murphy, S. A., Little,  R., Chakraborty, B. and Zhang, A.J. (2008). ``Screening experiments and the use of fractional factorial designs in behavioral intervention research'', \emph{American Journal of Public Health}, {\bf 98}, 1354$-$1359.

\item Nocedal, J. and Wright, S. J. (1999). {\em Numerical Optimization}. Springer, New York.

\item Pukelsheim, F. (1993). {\em Optimal Design of Experiments}, John Wiley \& Sons.

\item Pronzato, L. and Walter, E. (1988). ``Robust experiment design via maximin optimization'', {\it Math. Biosci.}, {\bf 89}, 161$-$176.

\item Rao, C. R. (1973). {\em Linear Statistical Inference and Its Applications}, John Wiley \& Sons, New York.

\item Russell, K.G., Woods, D.C., Lewis, S.M. and Eccleston, J.A. (2009). `` D-optimal designs for Poisson regression models'', \emph{Statistica Sinica}, {\bf 19}, 721$-$730.

%\item Seo, D. M., Goldschmidt-Clermont, P. J. and West, M. (2007). ``Of Mice And Men: Sparse Statistical Modeling in Cardiovascular Genomics", \emph{Annals of Applied Statistics}, {\bf 1}, 152$-$178.

\item Severin, V. (2000), ``Comparing statistical efficiency and respondent efficiency in choice experiments'', \emph{Ph.D. Thesis, University of Sydney}.

%\item Street, D. J. and Burgess, L. (2007). ``The Construction of Optimal Stated Choice Experiments'', \emph{Wiley Series in Probability and Statistics}.

\item Silvey, S. D., Titterington, D. M., and Torsney, B. (1978).  ``An algorithm for optimal designs on a finite design space", \emph{Commun. Stat. Theory Methods}, {\bf 14}, 1379$-$1389.

\item Stufken, J. and Yang, M. (2012). ``Optimal Designs for Generalized Linear Models'', In: \emph{Design and Analysis of Experiments}, Volume 3: Special Designs and Applications, K. Hinkelmann (ed.), Wiley, New~York.

%\item Stufken, J. and Yang, M. (2012). ``On locally optimal designs for generalized linear models with group effects'', \emph{Statistica Sinica}, {\bf 22}, 1765$-$1786.

\item Titterington, D. M.  (1976). ``Algorithms for computing D-optimal design on finite design spaces",  in \emph{Proc. of the 1976 Conf. on Information Science and Systems}, John~Hopkins University, {\bf 3}, 213$-$216.

\item Titterington, D. M. (1978). ``Estimation of correlation coefficients by ellipsoidal trimming", \emph{Appl. Stat.}, {\bf 27}, 227$-$234.

%\item Vasandani, V. and Govindaraj, T. (1995). ``Knowledge Organization in Intelligent Tutoring Systems for Diagnostic Problem Solving in Complex Dynamic Domains", \emph{IEEE Transactions on Systems, Man and Cybernetics}, 1076$-$1096.

\item Waterhouse, T. H., Woods, D. C., Eccleston, J. A. and Lewis, S. M. (2008). ``Design selection criteria for discrimination/estimation for nested models and a binomial response'', \emph{Journal of Statistical Planning and Inference}, {\bf 138}, 132$-$144.

\item Woods, D. C., Lewis, S. M., Eccleston, J. A. and Russell, K. G. (2006). ``Designs for generalized linear models with several variables and model uncertainty'', \emph{Technometrics}, {\bf 48}, 284$-$292.

\item Woods, D. C. and van de Ven, P. (2011). ``Blocked designs for experiments with non-normal response'', \emph{Technometrics}, {\bf 53}, 173$-$182.

\item Wynn, H.P. (1970). ``The sequential generation of D-optimum experimental designs'', {\em Annals of Mathematical Statistics}, {\bf 41}, 1655$-$1664.

\item Xu, H., Phoa, F. K. H. and Wong, W. K. (2009). ``Recent Developments in Nonregular Fractional Factorial Designs'', {\em Statistics Surveys}, {\bf 3}, 18$-$46.

\item Yang, J., Mandal, A. and Majumdar, D. (2012). ``Optimal designs for two-level factorial experiments with binary response'', \emph{Statistica Sinica}, {\bf 22}, 885$-$907.

\item Yang, M. and Stufken, J. (2009). ``Support points of locally optimal designs for nonlinear models with two parameters'', \emph{Annals of Statistics}, {\bf 37}, 518$-$541.

\item Yang, M., Zhang, B. and Huang, S. (2011). ``Optimal designs for generalized linear models with multiple design variables'', \emph{Statistica Sinica}, {\bf 21}, 1415$-$1430.

\item Yu, Y. (2010). ``Monotonic convergence of a general algorithm for computing optimal designs", \emph{Annals of Statistics}, {\bf 38}, 1593$-$1606.

\item Zangwill, W. (1969), ``Nonlinear Programming: A Unified Approach'', Prentice-Hall, New~Jersey.

\item Zayats, N. and Steinberg, D. M. (2010). ``Optimal Design of Experiments When Factors Affect Detection Capability", {\em Pakistan Journal of Statistics}, {\bf 26}, 15$-$37.
\end{description}

%\end{document}
\clearpage
\setcounter{page}{1}
\def\thepage{S\arabic{page}}

%\appendix

\fontsize{10.95}{14pt plus.8pt minus .6pt}\selectfont
\vspace{0.8pc}
\centerline{\large\bf OPTIMAL DESIGNS FOR $2^K$ FACTORIAL EXPERIMENTS}
\vspace{2pt}
\centerline{\large\bf WITH BINARY RESPONSE}
\vspace{.4cm}
\centerline{Jie Yang$^{1}$, Abhyuday Mandal$^{2}$ and Dibyen Majumdar$^{1}$}
\vspace{.4cm}
\centerline{\it  $^1$University of Illinois at Chicago and $^2$University of Georgia}
\vspace{.55cm}
\fontsize{9}{11.5pt plus.8pt minus .6pt}\selectfont

\setcounter{section}{1}
\setcounter{equation}{0}
\def\theequation{S\arabic{section}.\arabic{equation}}
\def\thesection{S\arabic{section}}

\renewcommand{\baselinestretch}{1.2}
\fontsize{10.95}{14pt plus.8pt minus .6pt}\selectfont

\begin{center} {\bf Supplementary Materials} \end{center}

\renewcommand{\theequation}{S.\arabic{equation}}

\noindent {\bf Connection between General Equivalence Theorem and Theorem~\ref{theorem30}:}

Extending the notations of this paper, we consider the problem when a design $\xi=\{({\bf x}_i, p_i),\ i=1, \ldots, 2^k\}$ maximizes the D-criterion $|M(\xi)|=|X'WX|$, where ${\bf x}_i$ is the $i$th row of $X$ and $X$ is the $2^k\times (d+1)$ {model matrix}.

\bigskip\noindent {\bf General Equivalence Theorem} (see, for example, Atkinson et.~al.~(2007)): $\xi$ maximizes $|M(\xi)|$ (or equivalently minimizes $\Psi\{M(\xi)\} = -\log|M(\xi)|$) if and only if
$$w_i{\bf x}_i'(X'WX)^{-1}{\bf x}_i \leq d+1$$ for each $i=1,\ldots, 2^k$ and equality holds if $p_i>0$.

Here's the outline of the proof of the General Equivalence Theorem described in Atkinson et.~al.~(2007, \S9.2, page 122): For each $i=1,\ldots, 2^k$, let $\bar{\xi}_i$ be the design supported only on ${\bf x}_i$, or in other words, it puts unit mass at the point ${\bf x}_i$ and let  $\xi'_i=(1-\alpha)\xi+\alpha\bar{\xi}_i$. The derivative of $\Psi$ in the direction $\bar{\xi}_i$ or ${\bf x}_i$ is $$\phi({\bf x}_i,\xi)=\lim_{\alpha \rightarrow 0^+} \frac{1}{\alpha} [\Psi(M(\xi'_i))-\Psi(M(\xi))]=(d+1)-w_i{\bf x}_i'(X'WX)^{-1}{\bf x}_i~.$$
Then $\xi$ is D-optimal if and only if $\min_i \phi({\bf x}_i,\xi)=0$ and $\phi({\bf x}_i, \xi)=0$ if $p_i>0$. Comparing with our proof of Theorem~\ref{theorem30}, $\xi'_i=(1-\alpha)\xi+\alpha\bar{\xi}_i = \xi + \alpha (\bar{\xi}_i-\xi)$ corresponds to our ${\bf p}_r + u\boldsymbol\delta_i^{(r)}$ with $u$ replaced by $\alpha$ and $\boldsymbol\delta_i^{(r)}$ replaced by $\bar{\xi}_i-\xi$. Therefore, $\phi({\bf x}_i, \xi)$ is equal to
$
\left.\frac{\partial f^{(r)}({\mathbf p}_r +u\boldsymbol\delta_i^{(r)})}{\partial u}\right|_{u=0}
$
and the if and only if condition comparing Atkinson et.~al.~(2007) becomes
\[
\left.\frac{\partial f^{(r)}({\mathbf p}_r +u\boldsymbol\delta_i^{(r)})}{\partial u}\right|_{u=0}
\begin{array}{cll}
=    & 0 & \mbox{ if }p_i >0\mbox{;}\\
\leq & 0 & \mbox{ otherwise.}
\end{array}
\]

The major difference between the general equivalence theorem and Theorem~3.11 is that the general equivalence theorem ends up with the inverse of $X'WX$, while we expressed the same set of conditions in terms of determinants with the aid of Lemma~3.1.1, as well as Lemma~S1.2 and Lemma~S1.3.
\hfill{$\Box$}

\bigskip\noindent {\bf Additional Results for Example 4.1:}
Consider a $2^3$ main-effects model with logit link. Suppose $\beta_1=0$. As a corollary of Theorem~\ref{theorem2^3},
the regular fractions $\{1,4,6,7\}$, $\{2,3,5,8\}$ are D-optimal half-fractions if and only
$$4\ \nu\left(|\beta_0| + |\beta_2| + |\beta_3|\right)
\geq \nu\left(|\beta_0| + |\beta_2| + |\beta_3| - 2\max_{0\leq i\leq 3} |\beta_i|\right).
$$ Note that $\nu(\eta) = \frac{1}{2 + e^\eta + e^{-\eta}}$ for logit link, which is symmetric about $0$.
To simplify the notations, let $\beta_{2\vee 3}=\max\{|\beta_2|,|\beta_3|\}$ and $\beta_{2\wedge 3}=\min\{|\beta_2|,|\beta_3|\}$.
The regular fractions $\{1,4,6,7\}$, $\{2,3,5,8\}$ are D-optimal half-fractions if and only
if one of three conditions below is satisfied:
\begin{eqnarray}
{\rm (i)} \mbox{\hspace{.2in}} |\beta_2|+|\beta_3|\leq \log 2; \label{eq:2^3b1=0}
\end{eqnarray}
\vspace{-.1in}
\begin{eqnarray*}
{\rm (ii)}  && |\beta_2|+|\beta_3| > \log 2,\>
       \beta_{2\vee 3} \leq \log\left(1+e^{-\beta_{2\wedge 3}}+\left[1+e^{-\beta_{2\wedge 3}}+e^{-2\beta_{2\wedge 3}}\right]^{1/2}\right),\nonumber\\
& & \mbox{and } |\beta_0| \leq \log\left(\frac{2\exp\{|\beta_2|+|\beta_3|\}-1}{\exp\{|\beta_2|+|\beta_3|\}-2}\right);\nonumber\\
{\rm (iii)}  && \beta_{2\vee 3} > \log\left(1+e^{-\beta_{2\wedge 3}}+\left[1+e^{-\beta_{2\wedge 3}}+e^{-2\beta_{2\wedge 3}}\right]^{1/2}\right),\nonumber \\
    &&   |\beta_{2\vee 3}| \leq \log\left(\frac{2e^{|\beta_{2\wedge 3}|}-1}{e^{|\beta_{2\wedge 3}|}-2}\right)
  \mbox{ and } |\beta_0| \leq \log\left(\frac{2e^{\beta_{2\vee 3}}-1}{e^{\beta_{2\vee 3}}-2}\right)-\beta_{2\wedge 3}.\nonumber
\end{eqnarray*}

 The above result is displayed in the right panel of Figure~\ref{fig0}. In the $x$- and $y$-axis, we have plotted $\beta_2$ and $\beta_3$ respectively. The rhomboidal region at the center (marked as $\infty$) represents the region where the regular fractions will always be D-optimal, irrespective of the values of $\beta_0$. The contours outside this region are for the upper bound of $|\beta_0|$.  Regular fractions will be D-optimal if the values of $|\beta_0|$ will be smaller than the upper bound with $\beta_2$ and $\beta_3$ falling inside the region outlined by the contour.

\section*{Proofs}

 We need two lemmas before the proof of Theorem~\ref{theorem30}.

\begin{lemma}\label{algo1lemma30}
Suppose ${\mathbf p}=(p_1,\ldots,p_{2^k})'$ satisfies $f\left({\mathbf p}\right)>0$. Given $i=1, \ldots, 2^k$,
\begin{eqnarray}\label{lem:f_i(x)}
f_i(z) = a_i z(1-z)^d+b_i (1-z)^{d+1},
\end{eqnarray}
for some constants $a_i$ and $b_i$.
If $p_i>0$, $b_i=f_i(0)$,
$a_i=\frac{f\left({\mathbf p}\right)-b_i\left(1-p_i\right)^{d+1}}{p_i\left(1-p_i\right)^d}$;
otherwise,
$b_i=f\left({\mathbf p}\right)$,
$a_i=f_i\left(\frac{1}{2}\right)\cdot 2^{d+1}-b_i$.
Note that $a_i\geq 0$, $b_i\geq 0$, and $a_i + b_i>0$.
\hfill{$\Box$}
\end{lemma}

\begin{lemma}\label{algo1lemma31}
Let $h(z)=az(1-z)^d+b(1-z)^{d+1}$ with $0\leq z \leq 1$ and $a\geq 0, b\geq 0, a+b >0$.
If $a>b(d+1)$, then
$\max_z h(z)=\left(\frac{d}{a-b}\right)^d\left(\frac{a}{d+1}\right)^{d+1}\mbox{ at }
z=\frac{a-b(d+1)}{(a-b)(d+1)}\>\><1.$
Otherwise, $\max_z h(z)=b$ at $z=0$.
\hfill{$\Box$}
\end{lemma}

\bigskip\noindent {\bf Proof of Theorem~\ref{theorem30}:} Note that $f({\mathbf p}) > 0$ implies $0\leq p_i <1$ for each $i=1,\ldots, 2^k$.
Since $\sum_i p_i = 1$, without any loss of generality, we assume $p_{2^k} >0$.
Define ${\mathbf p} _r=(p_1 ,\ldots,p_{2^k-1} )'$, and
$f^{(r)}({\mathbf p_r})=f(p_1,\ldots,p_{2^k-1},1-\sum_{i=1}^{2^k-1}p_i)$.

\smallskip
For $i=1, \ldots, 2^k-1$, let $\boldsymbol\delta_i^{(r)}=(-p _1,\ldots,-p _{i-1},1-p _i,-p _{i+1},\ldots,-p _{2^k-1})'$.
Then $f_i(z) =
f^{(r)}({\mathbf p} _r + u\boldsymbol\delta_i^{(r)})$ with $u=\frac{z-p_i }{1-p_i }$.
Since the determinant $|(\boldsymbol\delta_1^{(r)},\ldots,$ $\boldsymbol\delta_{2^k-1}^{(r)})| = p _{2^k} \neq 0$,
$\boldsymbol\delta_1^{(r)}, \ldots, \boldsymbol\delta_{2^k-1}^{(r)}$ are linearly independent
and thus may serve as a new basis of
\begin{equation}\label{eq:sr}
S_r = \{(p_1, \ldots, p_{2^k-1})'\ |\ \sum_{i=1}^{2^k-1} p_i \leq 1,\mbox{ and }
p_i \geq 0, i=1, \ldots, 2^k-1\}.
\end{equation}
Since $\log f^{(r)}({\mathbf p}_r)$ is concave, ${\mathbf p}_r$ maximizes
$f^{(r)}$ if and only if along each direction $\boldsymbol\delta_i^{(r)}$,
\[
\left.\frac{\partial f^{(r)}({\mathbf p}_r +u\boldsymbol\delta_i^{(r)})}{\partial u}\right|_{u=0}
=0\mbox{ if }p_i >0\mbox{;}\>\>\>
\leq 0\mbox{ otherwise.}
\]
That is, $f_i(z)$ attains its maximum at $z=p_i$, for each $i=1, \ldots, 2^k-1$ (and thus for $i=2^k$).
Based on Lemma~\ref{algo1lemma30} and Lemma~\ref{algo1lemma31}, it implies
one of the two cases:
\begin{itemize}
\item[(i)] $p_i=0$ and $f_i\left(\frac{1}{2}\right)\cdot 2^{d+1} - f({\mathbf p}) \leq f({\mathbf p}) (d+1)$;
\item[(ii)] $p_i > 0$, $a > b(d+1)$, and $a- b(d+1) = p_i (a-b)(d+1)$, where $b = f_i(0)$,
and $a= \frac{f({\mathbf p}) - b(1-p_i)^{d+1}}{p_i(1-p_i)^d}$.
\end{itemize}
The conclusion needed can be obtained by simplifying those two cases above.
\hfill{$\Box$}

\bigskip\noindent {\bf Proof of Theorem~\ref{dsaturatedconditiontheorem}:}
Let ${\mathbf p}_I$ be the minimally supported design satisfying $p_{i_1}=p_{i_2}=\cdots=p_{i_{d+1}}=\frac{1}{d+1}$.
Note that if $|X[i_1,i_2,\ldots,i_{d+1}]| = 0$, ${\mathbf p}_I$ can not be D-optimal.
Suppose $|X[i_1,i_2,\ldots,i_{d+1}]|\neq 0$, ${\mathbf p}_I$ is D-optimal if and only if ${\mathbf p}_I$ satisfies the conditions of Theorem~\ref{theorem30}.
By Lemma~\ref{xwx},
$f({\mathbf p}_I)= (d+1)^{-(d+1)} |X[i_1,i_2,\ldots,i_{d+1}]|^2 w_{i_1}w_{i_2}\cdots w_{i_{d+1}}$~.

 For $i\in {\mathbf I}$, $p_i=\frac{1}{d+1}$, $f_i(0)=0$. By case~(ii) of Theorem~\ref{theorem30},
$p_i = \frac{1}{d+1}$ maximizes $f_i(x)$.
For $i\notin {\mathbf I}$, $p_i=0$,
\begin{eqnarray*}
f_i\left(\frac{1}{2}\right) &=&[2(d+1)]^{-(d+1)} |X[i_1,\ldots,i_{d+1}]|^2 w_{i_1}\cdots w_{i_{d+1}}\\
                            &+&2^{-(d+1)}(d+1)^{-d} w_i\cdot w_{i_1}\cdots w_{i_{d+1}} \sum_{j\in {\mathbf I}} \frac{|X[\{i\}\cup {\mathbf I}\setminus\{j\}]|^2}{w_j}.
\end{eqnarray*}
Then $p_i=0$ maximizes $f_i(x)$ if and only if $f_i\left(\frac{1}{2}\right) \leq f({\mathbf p})
\frac{d+2}{2^{d+1}}$, which is equivalent to
$$
\sum_{j\in {\mathbf I}} \frac{|X[\{i\}\cup {\mathbf I}\setminus\{j\}]|^2}{w_j} \leq \frac{|X[i_1,i_2,\ldots,i_{d+1}]|^2}{w_i}.
$$
\hfill{$\Box$}

%%%%%%%%%%%%%%%%%%%%%%%%%%%%%%%%%%
%\iffalse

\bigskip\noindent {\bf Proof of Theorem~\ref{algo1theorem15}:}
Suppose the lift-one algorithm or its modified version converges at
${\mathbf p}^*=(p_1^*,\ldots,p_{2^k}^*)'$.
According to the algorithm, $|X'WX|>0$ at ${\mathbf p}^*$
and $p_i^*<1$ for $i=1,\ldots, 2^k$. The proof of Theorem~\ref{theorem30}
guarantees that ${\mathbf p}^*$ maximizes
$f({\mathbf p})=|X'WX|$.

%\smallskip
Now we show that the modified lift-one algorithm must converge to the maximum value $\max_{\mathbf p}|X'WX|$.
Based on the algorithm, we obtain a sequence of designs $\{{\mathbf p}_n\}_{n\geq 0} \subset
S_{r}$ defined in (\ref{eq:sr}) such that $|X'WX|>0$.
We only need to check the case when the sequence is infinite.
To simplify the notation, here we still denote $f({\mathbf p}) = f(p_1, \ldots, p_{2^k-1}, 1-\sum_{i=1}^{2^k-1} p_i)$
for ${\mathbf p} = (p_1, \ldots, p_{2^k-1})' \in S_r$.
Since that $f\left({\mathbf p}\right)$
is bounded from above on $S_{r}$ and $f\left({\mathbf p}_n\right)$
strictly increases with $n$, then $\lim_{n\rightarrow\infty} f\left({\mathbf p}_n\right)$ exists.

%\smallskip
Suppose $\lim_{n\rightarrow \infty}f\left({\mathbf p}_n\right) < \max_{\mathbf p}|X'WX|$.
Since $S_r$ is compact, there exists a ${\mathbf p_*}=(p_1^*,$ $\ldots,$ $p_{2^k-1}^*)'\in S_{r}$
and a subsequence $\{{\mathbf p}_{n_s}\}_{s\geq 1}
\subset \{{\mathbf p}_{10m}\}_{m\geq 0} \subset \{{\mathbf p}_{n}\}_{n\geq 0}$ such that
\[
0<f\left({\mathbf p}_*\right)
=\lim_{n\rightarrow \infty}f\left({\mathbf p}_n\right)
=\lim_{s\rightarrow \infty}f\left({\mathbf p}_{n_s}\right)
\mbox{ and }
\|{\mathbf p}_{n_s}-{\mathbf p}_*\|\longrightarrow 0\mbox{ as }
s\rightarrow\infty,
\]
where ``$\|\cdot\|$" represents the Euclidean distance.
Since ${\mathbf p}_*$ is not a solution maximizing $|X'WX|$,
by the proof of Theorem~\ref{theorem30} and the modified algorithm, there exists a $\boldsymbol\delta_i^{(r)}$ at ${\mathbf p}_*$
and an optimal $u_*\neq 0$ such that
${\mathbf p}_* + u_*\boldsymbol\delta_i^{(r)}\left({\mathbf p}_*\right) \in S_{r}$ and
$\Delta:= f\left({\mathbf p}_* + u_*\boldsymbol\delta_i^{(r)}\left({\mathbf p}_*\right)\right) - f\left({\mathbf p}_*\right) > 0$.

%\smallskip
As $s\rightarrow \infty$, ${\mathbf p}_{n_s}\rightarrow {\mathbf p}_*$,
its $i$th direction $\boldsymbol\delta_i^{(r)}\left({\mathbf p}_{n_s}\right)$ determined by the algorithm $\rightarrow \boldsymbol\delta_i^{(r)}\left({\mathbf p}_*\right)$,
and the optimal $u\left({\mathbf p}_{n_s}\right)$ $\rightarrow u_*$.
Thus
${\mathbf p}_{n_s} + u\left({\mathbf p}_{n_s}\right)\boldsymbol\delta_i^{(r)}\left({\mathbf p}_{n_s}\right)
\longrightarrow
{\mathbf p}_* + u_*\boldsymbol\delta_i^{(r)}\left({\mathbf p}_*\right)
$ and
\[
f\left({\mathbf p}_{n_s} + u\left({\mathbf p}_{n_s}\right)\boldsymbol\delta_i^{(r)}\left({\mathbf p}_{n_s}\right)\right)
- f\left({\mathbf p}_{n_s}\right)
\longrightarrow
f\left({\mathbf p}_* + u_*\boldsymbol\delta_i^{(r)}\left({\mathbf p}_*\right)\right)
- f\left({\mathbf p}_*\right)=\Delta.
\]
For all large enough $s$,
$f\left({\mathbf p}_{n_s} + u\left({\mathbf p}_{n_s}\right)\boldsymbol\delta_i^{(r)}\left({\mathbf p}_{n_s}\right)\right)
- f\left({\mathbf p}_{n_s}\right)
> \Delta/2 > 0
$.
However,
\[
f\left({\mathbf p}_{n_s} + u\left({\mathbf p}_{n_s}\right)\boldsymbol\delta_i^{(r)}\left({\mathbf p}_{n_s}\right)\right)
- f\left({\mathbf p}_{n_s}\right)
\leq f\left({\mathbf p}_{n_s+1}\right)
- f\left({\mathbf p}_{n_s}\right)
\leq
f\left({\mathbf p}_*\right)
- f\left({\mathbf p}_{n_s}\right)
\rightarrow 0
\]
The contradiction implies that
$\lim_{n\rightarrow \infty}f\left({\mathbf p}_n\right) = \max_{\mathbf p}|X'WX|$.
\hfill{$\Box$}
%\fi

\bigskip\noindent {\bf Proof of Theorem~\ref{theorem2^3}:}
Given $\beta_1=0$, we have
$w_1=w_5 = \nu(\beta_0+\beta_2+\beta_3)$,
$w_2=w_6 = \nu(\beta_0+\beta_2-\beta_3)$,
$w_3=w_7 = \nu(\beta_0-\beta_2+\beta_3)$,
$w_4=w_8 = \nu(\beta_0-\beta_2-\beta_3)$.
The goal is to find a half-fraction ${\mathbf I} = \{i_1, i_2, i_3, i_4\}$ which maximizes
$s({\mathbf I}) := |X[i_1, i_2, i_3, i_4]|^2w_{i_1} w_{i_2} w_{i_3} w_{i_4}$. For regular half-fractions ${\mathbf I} = \{1,4,6,7\}$
or $\{2,3,5,8\}$, $s({\mathbf I}) = 256 w_1 w_2 w_3 w_4$.
Note that $|X[i_1, i_2, i_3, i_4]|^2 = 0$ for $12$ half-fractions identified by $1=\pm A$, $1=\pm B$, $1=\pm C$, $1=\pm AB$, $1=\pm AC$, or $1=\pm BC$; and $|X[i_1, i_2, i_3, i_4]|^2=64$ for all other $56$ cases.

Without any loss of generality, suppose $w_1 \geq w_2 \geq w_3 \geq w_4$.
Note that the half-fraction $\{1, 5, 2, 6\}$ identified by $1=B$ leads to
$s({\mathbf I})=0$. Then the competitive half-fractions consist
of both $1$ and $5$, one element from the second block $\{2,6\}$, and one element from the third block $\{3,7\}$. The corresponding
$s({\mathbf I}) = 64 w_1^2w_2w_3$. In this case, the regular fractions are optimal ones if and only if $4w_4 \geq w_1$.
\hfill{$\Box$}

\bigskip We need the lemma below for Theorem~\ref{drobusttheorem}:
\begin{lemma}\label{indexsetlemma}
Suppose $k\geq 3$ and $\ d(d+1)\leq 2^{k+1}-4$. For any index set ${\mathbf I}=\{i_1,\ldots,i_{d+1}\}\subset \{1,\ldots,2^k\}$,
there exists another index set ${\mathbf I}'=\{i_1',\ldots,i_{d+1}'\}$
such that
\begin{equation}\label{equaldetcondition}
|X[i_1, \ldots, i_{d+1}]|^2 = |X[i_1', \ldots, i_{d+1}']|^2
\mbox{ and }
{\mathbf I}\cap {\mathbf I}'=\emptyset.
\end{equation}
\end{lemma}

\bigskip\noindent {\bf Proof of Lemma~\ref{indexsetlemma}:}
Note that $k\geq 3$ and $d(d+1)\leq 2^{k+1}-4$ imply
$d+1\leq 2^{k-1}$ and $\frac{d(d+1)}{2}<2^k-1$.
{
Let ${\mathbf I}=\{i_1,\ldots,i_{d+1}\}\subset \{1,\ldots,2^k\}$ be the given index set.
It can be verified that there exists a nonempty subset ${\mathbf J}\subset \{1,2,\ldots,k\}$, such that (i) the $i_1{\rm th}, \ldots, i_{d+1}{\rm th}$ rows of the matrix $[C_1, C_2,\ldots,$ $C_k]$ are same as the $i_1'{\rm th}, \ldots, i_{d+1}'{\rm th}$ rows of the matrix $[A_1, A_2, \ldots, A_k]$, where $A_1, \ldots, A_k$ are the columns of $X$ corresponding to the main effects, $C_i=-A_i$ if $i\in {\mathbf J}$ and $C_i=A_i$ otherwise; (ii) ${\mathbf I}'=\{i_1',\ldots,i_{d+1}'\}$ satisfies conditions (\ref{equaldetcondition}). Actually, the index set ${\mathbf I}'$ satisfying (i) always exists once ${\mathbf J}$ is given, since the $2^k$ rows of matrix $[A_1, \ldots, A_k]$ contain all possible vectors in $\{-1,1\}^k$. Then $|X[i_1, \ldots, i_{d+1}]|^2 = |X[i_1', \ldots,$ $i_{d+1}']|^2$ is guaranteed once ${\mathbf I}'$ satisfies (i). If ${\mathbf I}\cap {\mathbf I}'\neq \emptyset$, then there exists an $i_a' \in {\mathbf I}\cap {\mathbf I}'$ ($a \in \{1,\ldots, d+1\}$). Thus $i_a \in {\mathbf I}$ and the $i_a$th row of $[C_1, \ldots, C_k]$ is same as the $i_a'$th row of $[A_1, \ldots, A_k]$.  Based on the definitions of $C_1, \ldots, C_k$, the $i_a$th and $i_a'$th rows of $[A_1, \ldots, A_k]$ have the same entries at $A_i$ for all $i\notin {\mathbf J}$ but different entries at $A_i$ for all $i\in {\mathbf J}$. On the other hand, once the index pair $\{i_a, i_a'\} \subset {\mathbf I}$ is given, it uniquely determines the subset ${\mathbf J} \subset \{1, \ldots, k\}$.  Note that there are $2^k-1$ possible nonempty ${\mathbf J}$ but only $\frac{d(d+1)}{2}$ possible pairs in ${\mathbf I}$.  Since $\frac{d(d+1)}{2} < 2^k-1$, there is at least one ${\mathbf J}$ such that there is no pair in ${\mathbf I}$ corresponding to it.  For such a ${\mathbf J}$,  we must have ${\mathbf I}\cap {\mathbf I}'=\emptyset$ .
\hfill{$\Box$}
}

\bigskip\noindent {\bf Proof of Theorem~\ref{drobusttheorem}:}
Fixing any row index set $I=\{i_1,\ldots,i_{d+1}\}$ of $X$ such that
$|X[i_1,i_2,\ldots,i_{d+1}]|^2 > 0$,
among all the $(d+1)$-row fractional designs
satisfying $p_i=0, \forall i \notin {\mathbf I}$,
$|X'WX|$ attains its maximum
$\left(\frac{1}{d+1}\right)^{d+1}w_{i_1}\cdots w_{i_{d+1}}
\times |X[i_1,i_2,\ldots,i_{d+1}]|^2
$
at ${\mathbf p}_I$ satisfying $p_{i_1}=\cdots=p_{i_{d+1}}=\frac{1}{d+1}$. Given any other index set ${\mathbf I}'=\{i_1',\ldots,i_{d+1}'\}$ with minimally supported design ${\mathbf p_{I'}}$ satisfying $p_{i_1'}=\cdots=p_{i_{d+1}'}=\frac{1}{d+1}$, the loss of efficiency of ${\mathbf p}_{I}$ with respect to ${\mathbf p}_{I'}$ given ${\mathbf w}_{I'}=(w_1,\ldots,w_{2^k})'$ is
\begin{eqnarray*}
R_{I'}(I)
&=&1-\left(\frac{\psi({\mathbf p}_I,{\mathbf w}_{I'})}{\psi({\mathbf p}_{I'},{\mathbf w}_{I'})}\right)^{\frac{1}{d+1}}
=1-\left(\frac{w_{i_1}\cdots w_{i_{d+1}} |X[i_1,\ldots,i_{d+1}]|^2}{w_{i_1'}\cdots w_{i_{d+1}'} |X[i_1',\ldots,i_{d+1}']|^2}\right)^{\frac{1}{d+1}}\\
&\leq & 1-\frac{a}{b}\cdot \left(\frac{|X[i_1,\ldots,i_{d+1}]|^2}{ |X[i_1',\ldots,i_{d+1}']|^2}\right)^{\frac{1}{d+1}}.
\end{eqnarray*}
By Lemma~\ref{indexsetlemma}, there always exists
an index set ${\mathbf I}'=\{i_1',\ldots,i_{d+1}'\}$ such that
$|X[i_1',$ $\ldots,$ $i_{d+1}']|^2$  $=$  $|X[i_1,$ $\ldots,$ $i_{d+1}]|^2$
 and
${\mathbf I}\cap {\mathbf I}'=\emptyset$.
Let ${\mathbf w}_{I'}=(w_1,\ldots,w_{2^k})'$ satisfy $w_i = b, \forall i\in {\mathbf I}'$
and $w_i = a, \forall i\in {\mathbf I}$ (here we assume $(w_1, \ldots, w_{2^k})$ can take any point
in $[a,b]^{2^k}$).
Then the loss of efficiency of ${\mathbf p}_I$ with respect to this ${\mathbf w}_{I'}$ is at least $1-a/b$. If we choose ${\bf I}=\{i_1, \ldots, i_{d+1}\}$ which maximizes $X[i_1, \ldots, i_{d+1}]|^2$, then the corresponding ${\mathbf p}_I$ attains the minimum value $1-a/b$ of the maximum loss in efficiency compared to other minimally supported designs.
\hfill{$\Box$}

%\setcounter{section}{1}

%%%%%%%%%%%%%%%%%%%%%%%%%%%%%%%%%%

\bigskip\noindent We need two lemmas for the exchange algorithm for integer-valued allocations.

\begin{lemma}\label{algo3lemma1}
Let $g(z)=Az(m-z)+Bz+C(m-z)+D$ for real numbers $A > 0, B\geq 0, C\geq 0, D\geq 0$, and
integers $m>0, 0\leq z\leq m$.
Let $\Delta$ be the integer closest to $\frac{mA+B-C}{2A}$.
\begin{itemize}
\item[(i)] If $0\leq \Delta \leq m$, then
$\max_{0\leq z \leq m} g(z)=mC+D+(mA+B-C)\Delta-A\Delta^2
\mbox{ at }z=\Delta.$
\item[(ii)] If $\Delta < 0$, then $\max_{0\leq z\leq m}=mC+D$ at $z=0$.
\item[(iii)] If $\Delta > m$, then $\max_{0\leq z\leq m}=mB+D$ at $z=m$.
\end{itemize}
\end{lemma}

\begin{lemma}\label{algo3lemma2}
Let ${\mathbf n}=(n_1,\ldots,n_{2^k})'$, $W_n = {\rm diag}\{n_1w_1,\ldots,n_{2^k}w_{2^k}\}$,
$f({\mathbf n})=|X'W_nX|$. Fixing $1\leq i<j\leq 2^k$, let
\begin{eqnarray}\label{lem:67fij}
f_{ij}(z)&=&\nonumber
f\left(n_1,\ldots,n_{i-1},z,n_{i+1},\ldots,n_{j-1},m-z,n_{j+1},\ldots,n_{2^k}\right)\\
&\stackrel{\triangle}{=}& Az(m-z)+Bz+C(m-z)+D,
\end{eqnarray}
where $m=n_i+n_j$.
Then
{\it (i)} $D>0 \Longrightarrow B>0\mbox{ and }C>0$;
{\it (ii)} $B>0\mbox{ or }C>0 \Longrightarrow A>0$;
{\it (iii)} $f({\mathbf n})>0\Longrightarrow A>0$;
{\it (iv)} $D = f(n_1,\ldots,$ $n_{i-1},0,n_{i+1},\ldots,$ $n_{j-1},0,n_{j+1},\ldots,$ $n_{2^k})$.
{\it (v)} Suppose $m>0$, then
$A=\frac{2}{m^2}\left(2 f_{ij}\left(\frac{m}{2}\right)-f_{ij}(0)-f_{ij}(m)\right)$,
$B=\frac{1}{m}\left(f_{ij}(m)-D\right)$,
$C=\frac{1}{m}\left(f_{ij}(0)-D\right)$.
\end{lemma}

\subsection*{Exchange algorithm for real-valued allocations}

\begin{lemma}\label{algo2lemma2}
Let $g(z)=Az(e-z)+Bz+C(e-z)+D$ for nonnegative constants $A,B,C,D,e$.
Define $\Delta=\frac{eA+B-C}{2A}$.
\begin{itemize}
\item[(i)] If $0\leq \Delta \leq e$, then
$\max_{0\leq z \leq e} g(z)=eC+D+\frac{(eA+B-C)^2}{4A}
\mbox{ at }z=\Delta.$
\item[(ii)] If $\Delta < 0$, then $\max_{0\leq z\leq e}=eC+D$ at $z=0$.
\item[(iii)] If $\Delta > e$, then $\max_{0\leq z\leq e}=eB+D$ at $z=e$.
\end{itemize}
\end{lemma}

\begin{lemma}\label{algo2lemma3}
Let ${\mathbf p}=(p_1,\ldots,p_{2^k})'$,
$f({\mathbf p})=|X'WX|$, and
\begin{eqnarray*}f_{ij}(z) &:=&
f\left(p_1,\ldots,p_{i-1},z,p_{i+1},\ldots,p_{j-1},e-z,p_{j+1},\ldots,p_{2^k}\right)\\
&\stackrel{\triangle}{=}& Az(e-z)+Bz+C(e-z)+D,
\end{eqnarray*}
where $1\leq i<j\leq 2^k$ and $e=p_i+p_j$.
Then
{\it (i)} $D>0 \Longrightarrow B>0\mbox{ and }C>0$;
{\it (ii)} $B>0\mbox{ or }C>0 \Longrightarrow A>0$;
{\it (iii)} $f({\mathbf p})>0\Longrightarrow A>0$;
{\it (iv)} $D = f(p_1,\ldots,$ $p_{i-1},0,p_{i+1},\ldots,$ $p_{j-1},0,p_{j+1},\ldots,$ $p_{2^k})$;
{\it (v)} Suppose $e>0$, then
$A=\frac{2}{e^2}\left(2 f_{ij}\left(\frac{e}{2}\right)-f_{ij}(0)-f_{ij}(e)\right)$,
$B=\frac{1}{e}\left(f_{ij}(e)-D\right)$,
$C=\frac{1}{e}\left(f_{ij}(0)-D\right)$.
\end{lemma}

\bigskip\noindent {\bf Exchange algorithm for maximizing $f({\mathbf p})=f(p_1,\ldots,p_{2^k})=|X'WX|$}
\begin{itemize}
\item[$1^\circ$] Start with an arbitrary design ${\mathbf p}^{(0)}=(p^{(0)}_1,\ldots,p^{(0)}_{2^k})'$ such that $f({\mathbf p}^{(0)})>0$.
\item[$2^\circ$] Set up a random order of $(i,j)$ going through all pairs
$$\{(1,2),(1,3),\ldots,(1,2^k), (2,3), \ldots,(2^k-1,2^k)\}.$$
\item[$3^\circ$] For each $(i,j)$, if $e:=p_i^{(0)}+p_j^{(0)}=0$, let ${\mathbf p}^{(1)}={\mathbf p}^{(0)}$ and jump to $5^\circ$. Otherwise, let
\begin{eqnarray*}
f_{ij}(z)&=&f\left(p_1^{(0)},\ldots,p_{i-1}^{(0)},z,p_{i+1}^{(0)},\ldots,p_{j-1}^{(0)},e-z,p_{j+1}^{(0)},\ldots,p_{2^k}^{(0)}\right)\\
&=& A z(e-z) + Bz + C(e-z) + D
\end{eqnarray*}
with nonnegative constants $A,B,C,D$ determined by Lemma~\ref{algo2lemma3}.
\item[$4^\circ$] Define
${\mathbf p}^{(1)}
=\left(p_1^{(0)},\ldots,p_{i-1}^{(0)},z_*,p_{i+1}^{(0)},\ldots,p_{j-1}^{(0)},e-z_*,p_{j+1}^{(0)},\ldots,p_{2^k}^{(0)}\right)'
$
where $z_*$ maximizes $f_{ij}(z)$ with $0\leq z\leq e$ (see Lemma~\ref{algo2lemma2}). Note that
$f({\mathbf p}^{(1)})=f_{ij}(z_*)\geq f({\mathbf p}^{(0)})>0$.
\item[$5^\circ$] Repeat $2^\circ\sim 4^\circ$ until convergence (no more increase in terms of $f({\mathbf p})$
by any pairwise adjustment).
\end{itemize}

\begin{theorem}\label{algo2theorem3}
If the exchange algorithm converges, the converged ${\mathbf p}$
maximizes $|X'WX|$.
\end{theorem}

\bigskip\noindent {\bf Proof of Theorem~\ref{algo2theorem3}:} Suppose the exchange algorithm converges at
${\mathbf p}^*=(p_1^*,\ldots,p_{2^k}^*)'$.
According to the algorithm, $|X'WX|>0$ at ${\mathbf p}^*$.
Without any loss of generality, assume $p_{2^k}^*>0$.
Let ${\mathbf p}^*_r=(p_1^*,\ldots,p_{2^k-1}^*)$,
$l_r({\mathbf p}_r)=\log f_r({\mathbf p}_r)$, and
$f_r({\mathbf p}_r)
=f(p_1,$ $\ldots,$ $p_{2^k-1},$ $1-\sum_{i=1}^{2^k-1}p_i)$.
Then for $i=1,\ldots,2^k-1$,
$
\left.\frac{\partial l_r}{\partial p_i}\right|_{p^*_r}
=\frac{1}{f({\mathbf p}^*)}\cdot \left.\frac{\partial f_r}{\partial p_i}\right|_{p^*_r}
=0,\mbox{ if }p_i^*>0\mbox{;}\>\>\>
\leq 0,\mbox{ otherwise}
$.
Thus ${\mathbf p}^*$ (or ${\mathbf p}^*_r$) locally maximizes $l({\mathbf p})$
(or $l_r({\mathbf p}_r)$), and ${\mathbf p}^*$ attains
the global maximum of $f({\mathbf p})$ on $S$.
\hfill{$\Box$} \\

  Similar to the lift-one algorithm, we may
modify the exchange algorithm so that ${\mathbf p}^{(0)}$ won't be updated
until all potential pairwise exchanges among $p_i$'s have been checked. It can be verified that
the modified exchange algorithm must converge.

\end{document}